\documentclass[reqno,11pt]{article}
\usepackage{graphicx, amsmath, amsthm}
\usepackage{amssymb}

\input{stochbif.sty}

\begin{document}


\title{Metastability in Interacting Nonlinear \\
Stochastic Differential Equations~II: \\
Large-$N$ Behaviour}
\author{Nils Berglund, Bastien Fernandez and Barbara Gentz}
\date{}   

\maketitle

\begin{abstract}
\noindent
We consider the dynamics of a periodic chain of $N$ coupled overdamped
particles under the influence of noise, in the limit of large $N$. Each
particle is subjected to a bistable local potential, to a linear coupling
with its nearest neighbours, and to an independent source of white noise.
For strong coupling (of the order $N^2$), the system synchronises, in the
sense that all particles assume almost the same position in their
respective local potential most of the time. In a previous work, we showed
that the transition from strong to weak coupling involves a sequence of
symmetry-breaking bifurcations of the system's stationary configurations.
We analysed, for arbitrary $N$, the behaviour for coupling intensities
slightly below the synchronisation threshold. Here we describe the
behaviour for any positive coupling intensity $\gamma$ of order $N^2$,
provided the particle number $N$ is sufficiently large (as a function of
$\gamma/N^2$). In particular, we determine the transition time between
synchronised states, as well as the shape of the \lq\lq critical
droplet\rq\rq\, to leading order in $1/N$. Our techniques involve the
control of the exact number of periodic orbits of a near-integrable twist
map,  allowing us to give a detailed description of the system's potential
landscape, in which the metastable behaviour is encoded.  
\end{abstract}

\leftline{\small{\it Date.\/} November 21, 2006. Revised version, July 5, 2007.}
\leftline{\small 2000 {\it Mathematical Subject Classification.\/} 
37H20, 37L60 (primary), 37G40, 60K35 (secondary)}
\noindent{\small{\it Keywords and phrases.\/}
Spatially extended systems, 
lattice dynamical systems, 
open systems, 
stochastic differential equations, 
interacting diffusions, 
Ginzburg--Landau SPDE, 
transitions times,
most probable transition paths,
large deviations, 
Wentzell-Freidlin theory, 
diffusive coupling, 
synchronisation, 
meta\-stability, 
symmetry groups, 
symplectic twist maps.
}


\section{Introduction}
\label{sec_in}


In this paper, we continue our analysis of the metastable dynamics of a
periodic chain of coupled bistable elements, initiated in~\cite{BFG06a}. In
contrast with similar models involving discrete on-site variables, or
\lq\lq spins\rq\rq, whose metastable behaviour has been studied extensively
(see for instance~\cite{denHollander04,OlivieriVares05}), our model
involves continuous local variables, and is therefore described by a set of
interacting stochastic differential equations. 

The analysis of the metastable dynamics of such a system requires an
understanding of its $N$-dimensional \lq\lq potential landscape\rq\rq, in
particular the number and location of its local minima and saddles of index
$1$. In~\cite{BFG06a}, we showed that the number of stationary
configurations increases from $3$ to $3^N$ as the coupling intensity
$\gamma$ decreases from a critical value $\gamma_1$ of order $N^2$ to $0$.
This  transition from strong to weak coupling involves a sequence of
successive symmetry-breaking bifurcations, and we analysed in detail the
first of these bifurcations, which corresponds to desynchronisation. 

In the present work, we consider in more detail the behaviour for large
particle number $N$. In the limit $N\to\infty$, the system tends to a
Ginzburg--Landau stochastic partial differential equation (SPDE), studied
for instance in~\cite{EckmannHairer01,Rougemont02}. The Ginzburg--Landau
SPDE describes in particular the behaviour near bifurcation points of more
complicated equations, such as the stochastic Swift--Hohenberg
equation~\cite{BlomkerHairerPavliotis05}. For large but finite $N$, it
turns out that a technique known as \lq\lq spatial map\rq\rq\ analysis
allows us to obtain a precise control of the set of stationary points, for
values of the coupling well below the synchronisation threshold.  More
precisely, given a strictly positive coupling intensity $\gamma$ of order
$N^2$, there is an integer $N_0(\gamma/N^2)$ such that for all  $N\geqs
N_0(\gamma/N^2)$, we know precisely the number, location and type of the
potential's stationary points. This allows us to characterise the
transition times and paths between metastable states for all these values
of $\gamma$ and $N$. 

This paper is organised as follows. Section~\ref{sec_res} contains the
precise definition of our model, and the statement of all results. After
introducing the model in Section~\ref{ssec_mod} and describing general
properties of the potential landscape in Section~\ref{ssec_pot}, we explain
the heuristics for the limit $N\to\infty$ in Section~\ref{ssec_heuristics}.
In Section~\ref{ssec_largeN}, we state the detailed results on number and
location of stationary points for large but finite $N$, and in
Section~\ref{ssec_stoch} we present their consequences for the stochastic
dynamics. Section~\ref{sec_tm} contains the proofs of these results. The
proofs rely on a detailed analysis of the orbits of period $N$ of a
near-integrable twist map, which are in one-to-one correspondence with
stationary points of the potential. Appendix~\ref{app_ell} recalls some
properties of Jacobi's elliptic functions needed in the analysis, while
Appendix~\ref{sec_prtech} contains some more technical proofs of results
stated in Section~\ref{sec_tmgchi}. 

\subsection*{Acknowledgements}

Financial support by the French Ministry of Research, by way of the {\it
Action Concert\'ee Incitative (ACI) Jeunes Chercheurs, Mod\'elisation
stochastique de syst\`emes hors \'equilibre\/}, is gratefully acknowledged.
NB and BF thank the Weierstrass Institute for Applied Analysis and
Stochastics (WIAS), Berlin, for financial support and hospitality. BG
thanks the ESF Programme {\it Phase Transitions and Fluctuation Phenomena
for Random Dynamics in Spatially Extended Systems (RDSES)\/} for financial
support, and the Centre de Physique Th\'eorique (CPT), Marseille, for kind
hospitality. 


\section{Model and Results}
\label{sec_res}


\subsection{Definition of the Model}
\label{ssec_mod}

Our model of interacting bistable systems perturbed by noise is defined by
the following ingredients:
\begin{itemiz}
\item	The periodic one-dimensional lattice is given by $\lattice=\Z/N\Z$,
where $N\geqs2$ is the number of particles.

\item	To each site $i\in\lattice$, we attach a real variable $x_i\in\R$,
describing the position of the $i$th particle. The configuration space is
thus $\cX=\R^\lattice$.

\item	Each particle feels a local bistable potential, given by 
\begin{equation}
\label{mod1}
U(\xi) = \frac14 \xi^4 - \frac12 \xi^2\;,
\qquad
\xi\in\R\;.
\end{equation}
The local dynamics thus tends to push the particle towards one of the two
stable
positions $\xi=1$ or $\xi=-1$.

\item	Neighbouring particles in $\Lambda$ are coupled via a
discretised-Laplacian interaction, of intensity $\gamma/2$. 

\item	Each site is coupled to an independent source of noise, of
intensity $\sigma$.
The sources of noise are described by independent Brownian motions
$\set{B_i(t)}_{t\geqs0}$.
\end{itemiz}

The system is thus described by the following set of coupled stochastic
differential equations, defining a diffusion on $\cX$:
\begin{equation}
\label{mod2}
\6x^\sigma_i(t) = f(x^\sigma_i(t))\6t 
+ \frac\cng2 \bigbrak{x^\sigma_{i+1}(t)-2x^\sigma_i(t)+x^\sigma_{i-1}(t)}
\6t
+ \sigma \6B_i(t)\;,
\qquad i\in\Lambda\;,
\end{equation}
where the local nonlinear drift is given by 
\begin{equation}
\label{mod3}
f(\xi) = -\nabla U(\xi) = \xi - \xi^3\;.
\end{equation}
For $\sigma=0$, the system~\eqref{mod2} is a gradient system of the form
$\dot x = -\nabla V_\gamma(x)$, with potential 
\begin{equation}
\label{mod5}
V_{\gamma}(x) = 
\sum_{i\in\lattice} U(x_i) + \frac\cng4 \sum_{i\in\lattice}
(x_{i+1}-x_i)^2\;.
\end{equation}


\subsection{Potential Landscape and Metastability}
\label{ssec_pot}

The dynamics of the stochastic system depends essentially on the \lq\lq
potential landscape\rq\rq\ $V_\gamma$. 
As in~\cite{BFG06a}, we use the notations
\begin{equation}
\label{mod7}
\cS = \cS(\gamma)
= \setsuch{x\in\cX}{\nabla V_{\gamma}(x)=0}
\end{equation}
for the set of stationary points, and $S_k(\gamma)$ for the set of
$k$-saddles, that is, stationary points with $k$ unstable directions and
$N-k$ stable directions.

Understanding the dynamics for small noise essentially requires knowing the
graph $\cG=(\cS_0,\cE)$, in which two vertices $x^\star, y^\star\in\cS_0$
are connected by an edge $e\in\cE$ if and only if there is a $1$-saddle
$s\in\cS_1$ whose unstable manifolds converge to $x^\star$ and $y^\star$.
The system behaves essentially like a Markovian jump process on $\cG$. The
mean transition time from $x^\star$ to $y^\star$ is of order
$\e^{2H/\sigma^2}$, where $H$ is the potential difference between $x^\star$
and the lowest saddle leading to $y^\star$ (see~\cite{FW}).

It is easy to see that $\cS$ always contains at least the three points 
\begin{equation}
\label{pot1}
O = (0,\dots,0)\;, \qquad I^{\pm} = \pm(1,\dots,1)\;.
\end{equation}
Depending on the value of $\gamma$, the origin $O$ can be an $N$-saddle, or
a $k$-saddle for any odd $k$. The points $I^{\pm}$ always belong to
$\cS_0$, in fact we have
\begin{equation}
\label{pot2}
V_\gamma(x) > V_\gamma(I^+) = V_\gamma(I^-) = -\frac N4 \quad \forall
x\in\cX
\setminus\set{I^-,I^+}
\end{equation}
for all $\gamma>0$, so that $I^+$ and $I^-$ represent the most stable
configurations of the system. The three points $O$, $I^+$ and $I^-$ are the
only stationary points belonging to the diagonal
\begin{equation}
\label{pot2A}
\cD=\setsuch{x\in\cX}{x_1=x_2=\dots=x_N}\;. 
\end{equation}
On the other hand, being a
polynomial of degree $4$ in $N$ variables, the potential $V_\gamma$ can
have up to $3^N$ stationary points. 
 
The potential $V_\gamma(x)$, as well as the sets $S(\gamma)$ and
$S_k(\gamma)$, are invariant under the transformation group $G=G_N$ of
order $4N$ ($4$ if $N=2$), generated by the following three symmetries:
\begin{itemiz}
\item	the rotation around the diagonal given by
$R(x_1,\dots,x_N)=(x_2,\dots, x_N,x_1)$;
\item	the mirror symmetry $S(x_1,\dots,x_N)=(x_N,\dots,x_1)$;
\item	the point symmetry $C(x_1,\dots,x_N)=-(x_1,\dots,x_N)$.  
\end{itemiz}
\goodbreak

\begin{table}
\begin{center}
\begin{tabular}{|l|l|l|}
\hline
\vrule height 12pt depth 6pt width 0pt
$N$ & $x$ & Type of symmetry \\
\hline
\vrule height 12pt depth 6pt width 0pt
$4L$ & $A$ & 
$(x_1,\dots,x_L,x_L,\dots,x_1,-x_1,\dots,-x_L,-x_L,\dots,-x_1)$ \\
\vrule height 8pt depth 6pt width 0pt
     & $B$ & 
$(x_1,\dots,x_L,\dots,x_1,0,-x_1,\dots,-x_L,\dots,-x_1,0)$ \\
\hline 
\vrule height 12pt depth 6pt width 0pt
$4L+2$ & $A$ & 
$(x_1,\dots,x_{L+1},\dots,x_1,-x_1,\dots,-x_{L+1},\dots,-x_1)$ \\
\vrule height 8pt depth 6pt width 0pt
     & $B$ & 
$(x_1,\dots,x_L,x_L\dots,x_1,0,-x_1,\dots,-x_L,-x_L,\dots,-x_1,0)$ \\
\hline 
\vrule height 12pt depth 6pt width 0pt
$2L+1$ & $A$ & 
$(x_1,\dots,x_L,-x_L,\dots,-x_1,0)$ \\
\vrule height 8pt depth 6pt width 0pt
     & $B$ & 
$(x_1,\dots,x_L,x_L,\dots,x_1,x_0)$ \\
\hline 
\end{tabular} 
\end{center}
\caption[]
{Symmetries of the stationary points bifurcating from the origin at
$\gamma=\gamma_1$. The situation depends on whether $N$ is odd (in which
case we write $N=2L+1$) or even (in which case we write $N=4L$ or $N=4L+2$,
depending on the value of $N \pmod4$). Points
labelled $A$ are $1$-saddles near the desynchronisation bifurcation at
$\gamma=\gamma_1$, those labelled $B$ are $2$-saddles (for odd $N$, this is
actually a conjecture). More saddles of the same index are obtained by
applying elements of the symmetry group $G_N$ to $A$ and $B$.}
\label{table_desync}
\end{table}

In~\cite{BFG06a}, we proved the following results:
\begin{itemiz}
\item	There is a critical coupling intensity 
\begin{equation}
\label{pot2B}
\gamma_1 = \frac1{1-\cos(2\pi/N)}
\end{equation}
such that for all $\gamma\geqs\gamma_1$, the set of stationary points $\cS$
consists of the three points $O$ and $I^\pm$ only. The graph $\cG$ has two
vertices $I^\pm$, connected by a single edge. 

\item	As $\gamma$ decreases below $\gamma_1$, an even number of new
stationary points bifurcate from the origin. Half of them are $1$-saddles,
while the others are $2$-saddles. These points satisfy symmetries as shown
in Table~\ref{table_desync}. The potential difference
between $I^\pm$ and the $1$-saddles behaves like
$N(1/4-(\gamma_1-\gamma)^2/6)$ as $\gamma\nearrow\gamma_1$.

\item	New bifurcations of the origin occur for $\gamma=\gamma_M =
(1-\cos(2\pi M/N))^{-1}$, with $2\leqs M\leqs N/2$, in which saddles of
order higher than $2$ are created.
\end{itemiz}

The number of stationary points emerging from the origin at the
desynchronisation bifurcation at $\gamma=\gamma_1$ depends on the parity of
$N$. If $N$ is even, there are exactly $2N$ new points ($N$ saddles of
index $1$, and $N$ saddles of index $2$). If $N$ is odd, we were only able
to prove that the number of new stationary points is a multiple of $4N$,
but formulated the conjecture that there are exactly $4N$ stationary points
($2N$ saddles of index $1$, and $2N$ saddles of index $2$). We checked this
conjecture numerically for all $N$ up to $101$.  As we shall see in
Section~\ref{ssec_largeN}, the conjecture is also true for $N$ sufficiently
large. 


\subsection{Heuristics for the Large-$N$ Limit}
\label{ssec_heuristics}

We want to determine the structure of the set $\cS$ of stationary
points for large particle number $N$, and large coupling intensity
$\gamma$. For this purpose, we introduce the rescaled coupling intensity 
\begin{equation}
\label{LargeN1}
\gammat = \frac\gamma{\gamma_1}
= \frac{2\pi^2}{N^2} \gamma
\biggbrak{1+\biggOrder{\frac1{N^2}}}\;.
\end{equation}
Then, the desynchronisation bifurcation occurs for $\gammat=1$. 
We will consider values of $\gammat$ which may be smaller than $1$, but are
bounded away from zero. The reason why the set of stationary points can be
controlled in this regime is that as $N\to\infty$, the deterministic system
$\dot x=-\nabla V_\gamma(x)$ behaves like a Ginzburg--Landau partial
differential equation (PDE). Indeed, assume that the $N$ sites of the chain
are evenly distributed on a circle of radius $1$, and that there exists a
smooth function $u(\ph,t)$, $\ph\in\fS^1$, interpolating the coordinates of
$x(t)$ in such a way that 
\begin{equation}
\label{heur1}
u\Bigpar{2\pi\frac{i}{N},t} = x_i(t)
\qquad \forall i\in\Lambda\;.
\end{equation}  
Then in the limit $N\to\infty$, the
discrete Laplacian in~\eqref{mod2} converges to a constant times the second
derivative of $u(\cdot,t)$, and we obtain the PDE 
\begin{equation}
\label{heur2}
\sdpar ut(\ph,t) = f(u(\ph,t)) + \gammat \sdpar u{\ph\ph}(\ph,t)\;.
\end{equation}
Stationary solutions of~\eqref{heur2} satisfy the equation 
\begin{equation}
\label{heur3}
\gammat u''(\ph) = -f(u(\ph))\;,
\end{equation}  
describing the motion of a particle of mass $\gammat$ in the
\emph{inverted}\/ potential $-U(\ph)$. 
The prefactor $\gammat$ can be removed by scaling $\ph$: Setting
$u_0(\phi)=u(\sqrt{\gammat}\phi)$, we see that $u_0$ satisfies 
the equation
\begin{equation}
\label{heur4}
u_0'' = -f(u_0) = u_0^3 - u_0\;. 
\end{equation}
All periodic solutions of this equation are known
(cf.~Section~\ref{sec_tmaa}), and can be expressed in terms of Jacobi's
elliptic functions\footnote{For the reader's
convenience, we recall the definitions and main properties of Jacobi's
elliptic integrals and functions in Appendix~\ref{app_ell}.} as
\begin{equation}
\label{heur5}
u_0(\phi) = a(\kappa) \sn
\biggpar{\frac{\phi-\phi_0}{\sqrt{1+\kappa^2}},\kappa}\;,
\end{equation} 
where
\begin{itemiz}
\item 	$\phi_0$ is an arbitrary phase;
\item	$\kappa\in[0,1)$ is an auxiliary parameter controlling the shape of
the function: For small $\kappa$, the function is close to a sine, while it
approaches a square wave as $\kappa\nearrow1$;
\item	the amplitude $a(\kappa)$ is given by 
\begin{equation}
\label{heur6}
a(\kappa)^2 = \frac{2\kappa^2}{1+\kappa^2}\;;
\end{equation} 
\item	the period of $u_0(\ph)$ is $4\sqrt{1+\kappa^2}\JK(\kappa)$, where
$\JK$ denotes the complete elliptic integral of the first kind.
\end{itemiz}

\begin{figure}
\centerline{
\includegraphics*[clip=true,height=45mm]{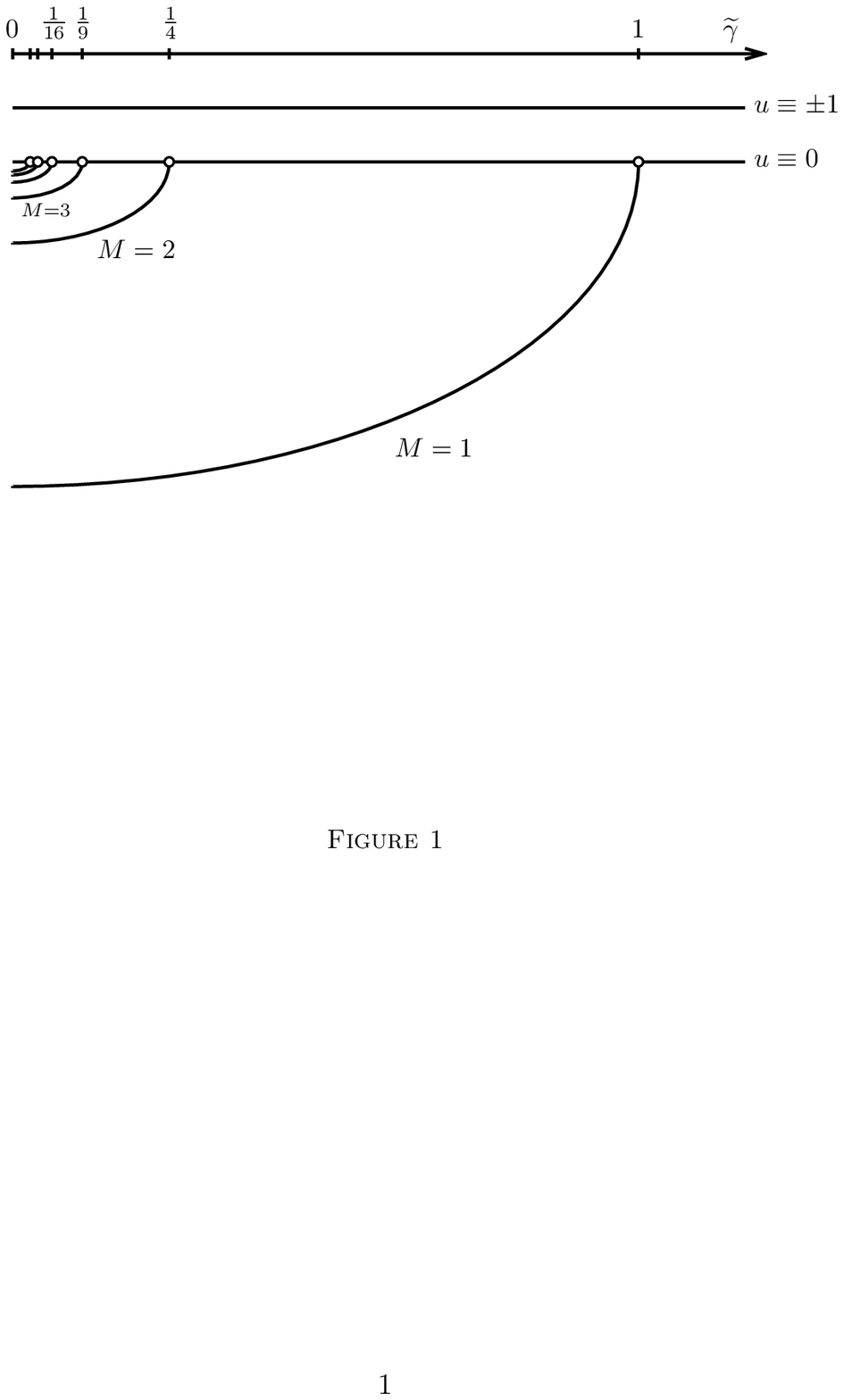}
\hspace{3mm}
\includegraphics*[clip=true,height=45mm]{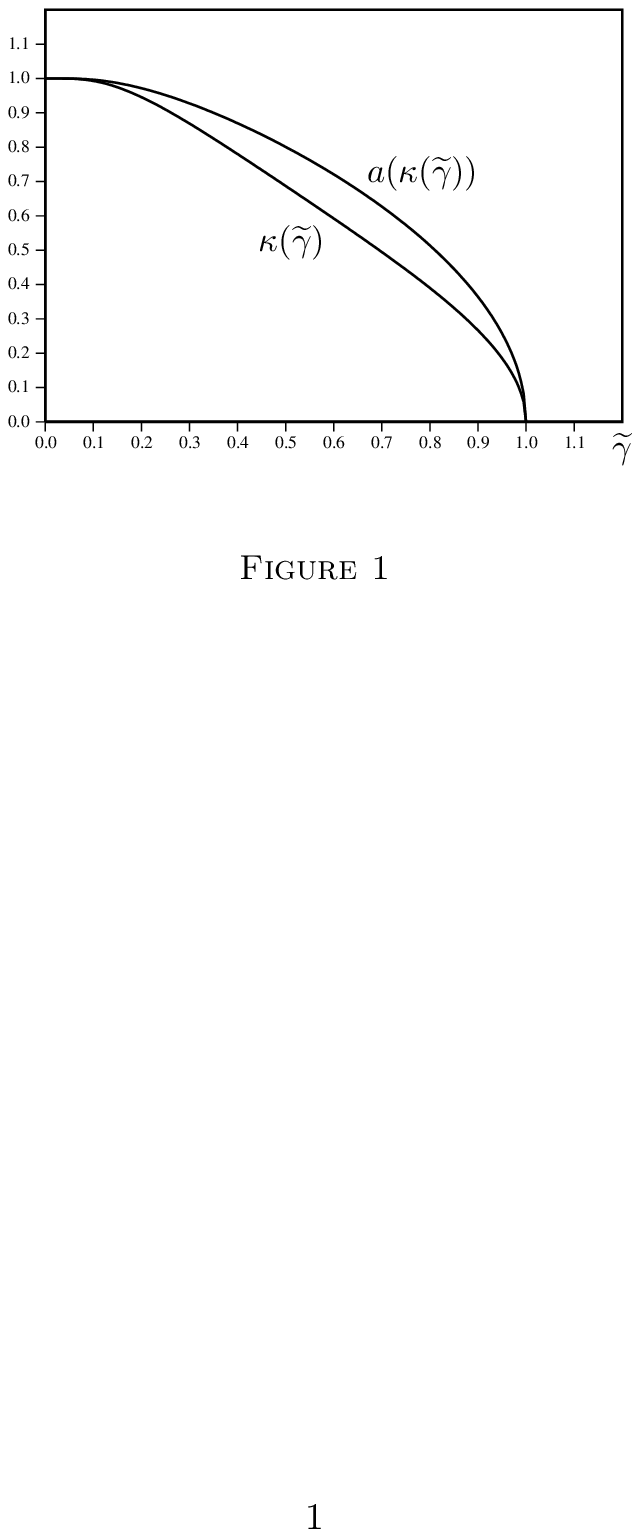}
}
 \figtext{ 
 	\writefig	-0.15	4.7	{\footnotesize{{\bf (a)}}}
  	\writefig	7.5	4.7	{\footnotesize{{\bf (b)}}}
}
 \caption[]
 {{\bf (a)}~Schematic bifurcation diagram of the limiting
PDE~\eqref{heur3}. Whenever the rescaled coupling intensity $\gammat$
decreases below $1/M^2$, $M=1,2,\dots$, a one-parameter family of
stationary solutions bifurcates from the identically zero solution
$u\equiv0$. 
{\bf (b)}~Relations between the shape parameter $\alkappa$, the
amplitude $a$ and
the rescaled coupling intensity $\gammat$ for winding number $M=1$.}
\label{fig_largeN1}
\end{figure}

We are looking for periodic solutions of~\eqref{heur4} of period
$2\pi/\sqrt{\gammat}$. Such solutions exist whenever the shape 
parameter $\kappa$ satisfies the condition
\begin{equation}
\label{heur7}
4\sqrt{1+\kappa^2}\JK(\kappa) = \frac{2\pi M}{\sqrt{\gammat}}
\end{equation} 
for some integer $M$, which plays the r\^ole of a \emph{winding number}\/
controlling the number of sign changes of $u_0$. Equation~\eqref{heur7}
imposes a relation between shape parameter $\kappa$ and rescaled coupling
intensity $\gammat$, shown in \figref{fig_largeN1} in the case $M=1$. On
the other hand, the phase $\phi_0$ is completely free.

The left-hand side of~\eqref{heur7} being bounded below by $2\pi$,
solutions of given winding number $M$ exist provided $\gammat \leqs 1/M^2$.
The smaller $\gammat$, the more different types of periodic solutions
exist. A new one-parameter family of stationary solutions, parametrised by
$\phi_0$, bifurcates from the identically zero solution every time
$\gammat$ becomes smaller than $1/M^2$, $M=1,2,\dots$
(\figref{fig_largeN1}a). 

Finally, note that for stationary points $x$ satisfying~\eqref{heur1}, with
$u$ given by~\eqref{heur5}, the value of the renormalised potential
$V_\gamma(x)/N$ converges, in the limit $N\to\infty$, to an integral which
can be computed explicitly (see Section~\ref{sec_tmaa}) in terms of the
parameter $\kappa$:
\begin{equation}
\label{heur8}
\lim_{N\to\infty} \frac{V_\gamma(x)}{N} 
= -\frac1{3(1+\alkappa^2)} \biggbrak{\frac{2+\alkappa^2}{1+\alkappa^2} - 
2\frac{\JE(\alkappa)}{\JK(\alkappa)}}\;, 
\end{equation} 
where $\JE$ denotes the complete elliptic integral of the second kind. 

If we were to add noise to the PDE~\eqref{heur2}, we would obtain a
Ginzburg--Landau SPDE. In that case, we expect that the
configurations of highest energy reached in the course of a typical
transition from $u\equiv-1$ to $u\equiv1$ are of the form~\eqref{heur5},
with winding number $M=1$. As a consequence, the potential difference
in~\eqref{heur8} should be governing the typical time of such transitions.
However, proving this would involve an infinite-dimensional version of
Wentzell--Freidlin theory, moreover in a degenerate situation, which is
beyond the scope of the present work. We will henceforth consider
situations with large, but finite particle number. 


\subsection{Main Results: Stationary Points for Large but Finite $N$}
\label{ssec_largeN}

We examine now the structure of the set $\cS(\gamma)$ for large, but finite
particle number $N$. Instead of the limiting differential
equation~\eqref{heur3}, stationary points satisfy the difference equation 
\begin{equation}
\label{largeN01}
\frac\cng2 \bigbrak{x_{n+1}-2x_n+x_{n-1}} = -f(x_n)\;, 
\qquad n\in\Lambda\;.
\end{equation}
The key idea of the analysis is to interpret $n$ as discrete time, and to
consider~\eqref{largeN01} as defining $x_{n+1}$ in terms of $x_n$ and
$x_{n-1}$. Setting $v_n = x_n - x_{n-1}$ allows to rewrite~\eqref{largeN01}
as the system 
\begin{equation}
\label{largeN02}
\begin{split}
x_{n+1} &= x_n + v_{n+1}\;, \\
v_{n+1} &= v_n - 2\gamma^{-1} f(x_n)\;. 
\end{split}
\end{equation}
The map $(x_n,v_n)\mapsto(x_{n+1},v_{n+1})$ is an area-preserving twist map
(\lq\lq twist\rq\rq\ meaning that $x_{n+1}$ is a monotonous  function
of~$v_n$), for the study of which many tools are available~\cite{Meiss92}.
Stationary points of the potential $V_\gamma$ are in one-to-one
correspondence with periodic points of period $N$ of this map. If we
further scale $v$ by a factor $\eps=\sqrt{2/\gamma}$, we obtain the
equivalent map  
\begin{equation}
\label{largeN03}
\begin{split}
x_{n+1} &= x_n + \eps y_{n+1}\;, \\
y_{n+1} &= y_n - \eps f(x_n)\;. 
\end{split}
\end{equation}
The regime of large particle number $N$ and finite rescaled coupling
intensity $\gammat$ corresponds to large $\gamma$, and thus to small
$\eps$. The map~\eqref{largeN03} is a discretisation of the system of
ordinary differential equations $\dot x=y$, $\dot y=-f(x)$, which is
equivalent to the continuous limit equation~\eqref{heur4}. There
should thus be some similarity between the orbits of the
map~\eqref{largeN03} and of the system~\eqref{heur4}. In particular, one
easily checks that the energy 
\begin{equation}
\label{largeN04}
E(x,\dot x) = \frac12\dot x^2 - U(x)\;,
\end{equation}
which is conserved in the continuous limit, changes only slightly, by an
amount of order $\eps^2$, for the map~\eqref{largeN03} (setting $y=\dot
x$). The map is thus close to integrable, which makes its analysis
accessible to perturbation theory. 

\begin{figure}
\centerline{\includegraphics*[clip=true,width=142mm]{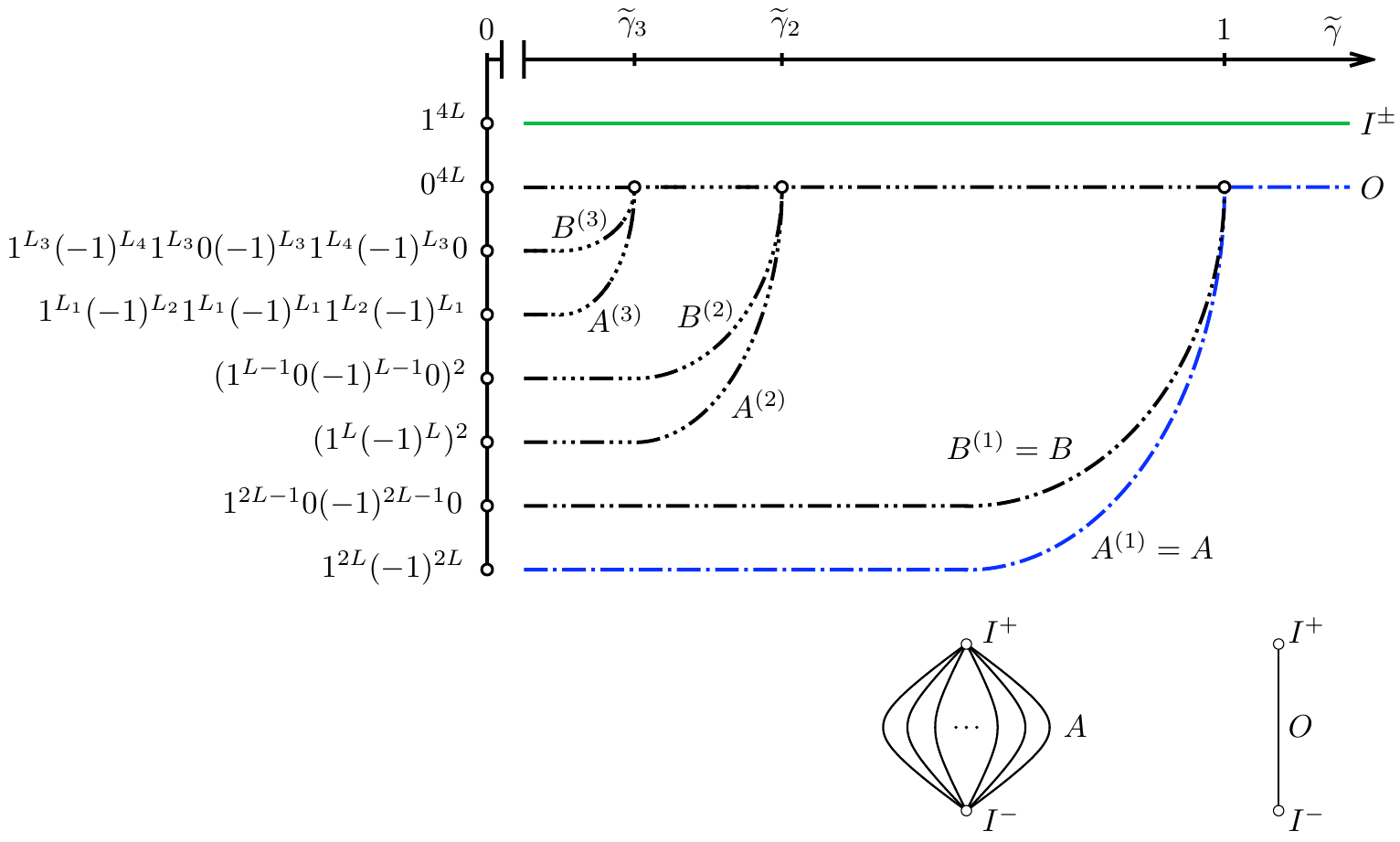}}
 \figtext{ }
 \caption[]
 {Partial bifurcation diagram for a case where $N=4L$ is a multiple of
 four,  and some associated graphs $\cG$. Only one stationary point per
 orbit of the symmetry group $G$ is shown. Dash--dotted curves with $k$
 dots
 represent $k$-saddles.  The symbols at the left indicate the zero-coupling
 limit of the stationary points' coordinates, for instance
 $1^{2L}(-1)^{2L}$ stands for a point whose first $2L$ coordinates are
 equal to $1$, and whose last $2L$ coordinates are equal to $-1$. The
 numbers associated with the branch created at $\gammat_3$ are
 $L_1=\intpart{2L/3}$, $L_2=2(L-L_1)$,  $L_3=\intpart{2L/3+1/2}$ and
 $L_4=2(L-L_3)-1$ (in case $N$ is a multiple of $12$, there are more
 vanishing coordinates).}
\label{fig_bifN}
\end{figure}

Our main result, obtained by analysing the map~\eqref{largeN03}, is
that the bifurcation diagram looks like the one shown
in~\figref{fig_bifN}. Namely,
\begin{itemiz}
\item 	For $\gammat > 1$, $I^+$, $I^-$ and $O$ are the only stationary
points.
\item	Below $\gammat = 1$, an explicitly known number of saddles of
index $1$ and $2$ bifurcate from the origin. 
\item	For any $2\leqs M\leqs N/2$, an explicitly known number of saddles
of index $2M-1$ and $2M$ bifurcate from the origin at $\gammat =
\gammat_M$, where 
\begin{equation}
\label{LargeN2}
\gammat_M = \frac{\gamma_M}{\gamma_1} 
= \frac{1-\cos(2\pi/N)}{1-\cos(2\pi M/N)} 
= \frac1{M^2} + \biggOrder{\frac1{N^2}}\;.
\end{equation}
\item	For any fixed $M$, if $N$ is sufficiently large compared to $M$,
the above list of stationary points is complete for $\gammat>\gammat_M$. In
particular, there are no secondary bifurcations of existing branches of
stationary points, and no stationary points created by saddle--node
bifurcation for these values of $\gammat$. 
\end{itemiz}

The main difficulty is to rule out the appearance of other stationary
points away from the origin. Indeed, for perturbed integrable maps it is
easy to obtain a lower bound on the number of periodic points, using the
Poincar\'e--Birkhoff theorem, but it is hard to obtain an upper bound. One
might imagine a scenario where stationary points appear far from the
origin, which ultimately offer a more economic path for the transition
from $I^-$ to $I^+$. 

We now give the precise formulation of the results. We first describe
the behaviour between the first two bifurcation values $\gammat_1$ and
$\gammat_2$. Below, $\gcd(a,b)$ denotes the greatest common divisor of
two integers $a$ and $b$, and $O_x=\setsuch{gx}{g\in G}$ denotes the group
orbit of a point $x\in\cX$ under the symmetry group $G$.

\begin{theorem}
\label{thm_LargeN1}
There exists $N_1<\infty$ such that when $N\geqs N_1$ and
$\gammat_2<\gammat<\gammat_1=1$, the set $\cS$ of stationary points of the
potential $V_\gamma$ has cardinality 
\begin{equation}
\label{LargeN3B}
\abs{\cS} = 3 + \frac{4N}{\gcd(N,2)} = 
\begin{cases}
3+2N & \text{if $N$ is even\;,}\\
3+4N & \text{if $N$ is odd\;,}
\end{cases}
\end{equation}
There exist points $A=A(\gammat)$ and $B=B(\gammat)$ in $\cX$ such that
$\cS$ can be decomposed as\footnote{If $N$ is even, the orbits $O_A$ and
$O_B$ contain $N$ instead of $2N$ points, because $R^{N/2}=-\one$.} 
\begin{align}
\nonumber
\cS_0 &= O_{I^+} = \set{I^+,I^-}\;, \\
\nonumber
\cS_1 &= O_{A} = \set{\pm A, \pm RA, \dots, \pm R^{N-1}A}\;, \\
\nonumber
\cS_2 &= O_{B} = \set{\pm B, \pm RB, \dots, \pm R^{N-1}B}\;, \\
\cS_3 &= O_{O} = \set{O}\;.
\label{LargeN4}
\end{align} 
The potential difference between the $1$-saddles and the well bottoms 
(which is the same for all $1$-saddles and well bottoms) satisfies
\begin{equation}
\label{largeN6}
\frac{V_\gamma(A(\gammat))-V_\gamma(I^\pm)}N = \frac14 -
\frac1{3(1+\alkappa^2)}
\biggbrak{\frac{2+\alkappa^2}{1+\alkappa^2}  - 2 
\frac{\JE(\alkappa)}{\JK(\alkappa)}} + 
\biggOrder{\frac{\alkappa^2}{N}}\;,
\end{equation}
where $\alkappa=\alkappa(\gammat)$ is defined implicitly by 
\begin{equation}
\label{LargeN3}
\gammat = \frac{\pi^2}{4\JK(\alkappa)^2(1+\alkappa^2)}\;.
\end{equation}

\end{theorem}

The detailed proofs are given in Section~\ref{sec_tm}. 

The second result, which is also proved in Section~\ref{sec_tm}, concerns
the behaviour for subsequent bifurcation values $\gammat_M$, $M\geqs 2$. 

\begin{theorem}
\label{thm_LargeN2}
For any $M\geqs2$,  there exists $N_M<\infty$ such that when $N\geqs N_M$
and $\gammat_{M+1}<\gammat<\gammat_M$, the set $\cS$ of stationary points
of the potential $V_\gamma$ has cardinality 
\begin{equation}
\label{LargeN7}
\abs{\cS} = 3 + \sum_{m=1}^M \frac{4N}{\gcd(N,2m)}\;.
\end{equation}
There exist points $A^{(m)}$ and $B^{(m)}$ in $\cX$, $m=1,\dots,M$, such
that $\cS$ can be decomposed as 
\begin{align}
\nonumber
\cS_0 &= O_{I^+} = \set{I^+,I^-}\;, \\
\nonumber
\cS_{2m-1} &= O_{A^{(m)}}\;, 
&
m &= 1,\dots,M\;, \\
\nonumber
\cS_{2m} &= O_{B^{(m)}}\;, 
&
m &= 1,\dots,M\;, \\
\cS_{2M+1} &= O_{O} = \set{O}\;,
\label{LargeN8}
\end{align}  
The potential difference between the saddles $A^{(m)}_j(\gammat)$ and the
well bottoms  satisfies
a similar relation as~\eqref{largeN6}, but with 
$\alkappa=\alkappa(m^2\gammat)$.
\end{theorem}

\begin{figure}
\centerline{\includegraphics*[clip=true,height=40mm]{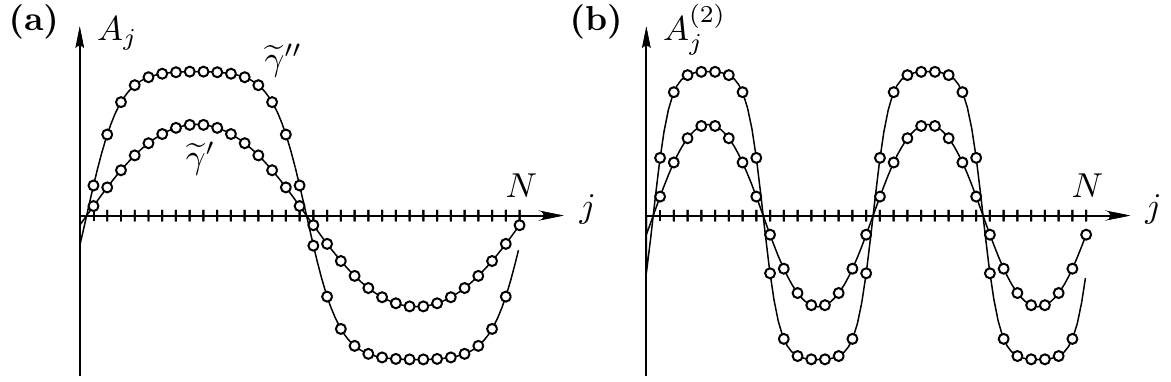}}
 \figtext{ }
 \caption[]
 {{\bf (a)} Coordinates of the $1$-saddles $A$ in the case 
 $N=32$, shown for two different values of the coupling
$\gammat'>\gammat''$.
 {\bf (b)} Coordinates of the $3$-saddles $A^{(2)}$ in the case 
 $N=32$, shown for the coupling intensities $\gammat'/4$ and
$\gammat''/4$.}
\label{fig_symell}
\end{figure}

\begin{figure}
\centerline{\includegraphics*[clip=true,width=120mm]{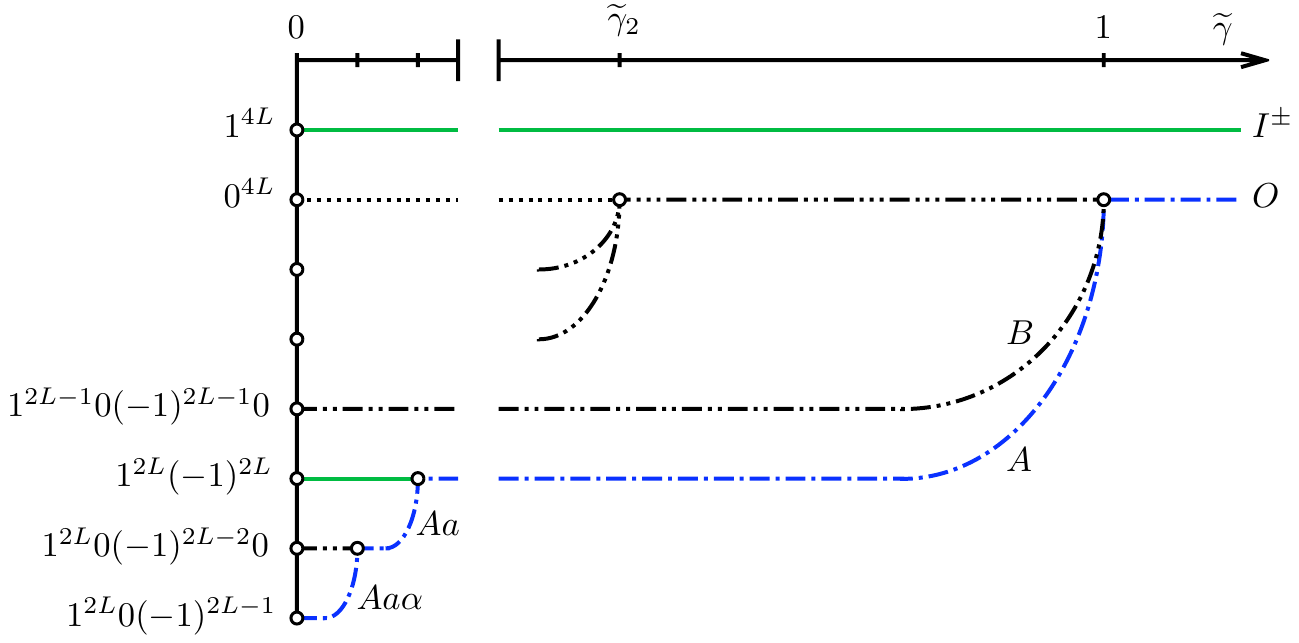}}
 \figtext{ }
 \caption[]
 {Partial bifurcation diagram for a case where $N=4L$ is a multiple of
four,
 showing the expected bifurcation behaviour of the critical $1$-saddle 
 in the zero-coupling limit.}
\label{fig_bifN2}
\end{figure}

\begin{remark}
The proof actually yields information on the coordinates of 
the points $A=A(\gammat)$ and $B=B(\gammat)$:
\begin{itemiz}
\item	The coordinates of $A$ and $B$ satisfy the symmetries indicated
in Table~\ref{table_desync}.

\item	If $N$ is even, the coordinates of $A$ and $B$ are given by 
\begin{align}
\nonumber
A_j(\gammat) &= a(\alkappa(\gammat)) \sn
\biggpar{\frac{4\JK(\alkappa(\gammat))}N \bigpar{j-\tfrac12},
\alkappa(\gammat)}  +
\biggOrder{\frac1{N}}\;, \\
B_j(\gammat) &= a(\alkappa(\gammat)) \sn
\biggpar{\frac{4\JK(\alkappa(\gammat))}N j, \alkappa(\gammat)}  +
\biggOrder{\frac1{N}}\;,
\label{LargeN5}
\end{align}
where the amplitude $a(\alkappa)$ is the one defined in~\eqref{heur6}. 

\item	If $N$ is odd, the coordinates of $A$ and $B$ are given by
\begin{align}
\nonumber
A_j(\gammat) &= a(\alkappa(\gammat)) \sn
\biggpar{\frac{4\JK(\alkappa(\gammat))}N j, \alkappa(\gammat)}  +
\biggOrder{\frac1{N}}\;, \\
B_j(\gammat) &= a(\alkappa(\gammat)) \cn
\biggpar{\frac{4\JK(\alkappa(\gammat))}N j, \alkappa(\gammat)}  +
\biggOrder{\frac1{N}}\;.
\label{largeN7}
\end{align}

\item	The components of $A^{(m)}$ and $B^{(m)}$ are given by similar
expressions, with $\gammat$ replaced by $m^2\gammat$, $j-\frac12$ replaced
by $m(j-\frac12)$ and $j$ replaced by $mj$.

\item	Note that the total number of stationary points accounted for by
these results is of the order $N^2$, which is much less than the $3^N$
points present at zero coupling. Many additional stationary points thus
have to be created as the rescaled coupling intensity $\cngt$ decreases
sufficiently, either by pitchfork-type second-order bifurcations of already
existing points, or by saddle-node bifurcations. However, the values
$\gammat(N)$ for which these bifurcations occur have to satisfy
$\lim_{N\to\infty}\gammat(N)=0$. 

The existence of second-order bifurcations
follows from stability arguments. For instance, for even $N$, the point
$A(\gammat)$ converges to $(1,1,\dots,1,-1,-1,\dots,-1)$ as $\gammat\to0$,
which is a local minimum of $V_\gamma$ instead of a $1$-saddle. The
$A$-branch thus has to bifurcate at least once as the coupling decreases to
zero (\figref{fig_bifN2}). For odd $N$, by contrast, the point $A(\gammat)$
converges to $(1,1,\dots,1,0,-1,-1,\dots,-1)$ as $\gammat\to0$, which is
also a $1$-saddle. We thus expect that the point $A(\gammat)$ does not
undergo any bifurcations for $0\leqs\gammat<1$ if $N$ is odd. 
\end{itemiz} 
\end{remark}

\subsection{Stochastic Case}
\label{ssec_stoch}

We return now to the behaviour of the system of stochastic differential
equations
\begin{equation}
\label{stoch1}
\6x^\sigma_i(t) = f(x^\sigma_i(t))\6t 
+ \frac\cng2 \bigbrak{x^\sigma_{i+1}(t)-2x^\sigma_i(t)+x^\sigma_{i-1}(t)}
\6t
+ \sigma \6B_i(t)\;, 
\qquad i\in\Lambda\;.
\end{equation}
Our main goal is to characterise the noise-induced transition from the
configuration $I^-=(-1,-1,\dots,-1)$ to the configuration
$I^+=(1,1,\dots,1)$. In particular, we are interested in the time needed
for this transition to occur, and by the shape of the critical 
configuration, i.e., the configuration of highest energy reached during the
transition. 

In~\cite[Theorem~2.7]{BFG06a}, we obtained that in the synchronisation
regime $\gammat>1$, for any initial condition $x_0$ in a ball $\cB(I^-,r)$
of radius $r<1/2$ around $I^-$, any particle number $N\geqs2$ and any
constant $\delta>0$, the first-hitting time $\tau_+ = \tauhit(\cB(I_+,r))$
of a ball $\cB(I^+,r)$ of radius $r$ around $I^+$ satisfies 
\begin{equation}
\label{stoch4}
\lim_{\sigma\to0}
\bigprobin{x_0}{\e^{(N/2-\delta)/\sigma^2} < \tau_+ <
\e^{(N/2+\delta)/\sigma^2}} = 1
\end{equation}
and 
\begin{equation}
\label{stoch4b}
\lim_{\sigma\to0} \sigma^2 \log 
\expecin{x_0}{\tau_+} = \frac N2\;.
\end{equation}
This means that in the synchronisation regime, the transition between $I^-$
and $I^+$ takes a time of the order $\e^{N/2\sigma^2}$. Furthermore, for
any fixed radius $R\in(r,1/2)$, the first-hitting time $\tau_O =
\tauhit(\cB(O,r))$ of a ball around the origin satisfies 
\begin{equation}
\label{stoch6}
\lim_{\sigma\to0}
\bigpcondin{x_0}{\tau_O < \tau_+}{\tau_+ < \tau_-} = 1\;,
\end{equation}
where $\tau_- = \inf\setsuch{t>\tauexit(\cB(I^-,R))}{x_t\in\cB(I^-,r)}$ is
the time of first return to the small ball $\cB(I^-,r)$ after leaving the
larger ball $\cB(I^-,R)$. This means that during a transition, the system
is likely to pass close to the origin, i.e., the origin, being the only
saddle of $V_\gamma$, is the critical configuration of the transition.
 
\begin{figure}
\centerline{\includegraphics*[clip=true,height=55mm]{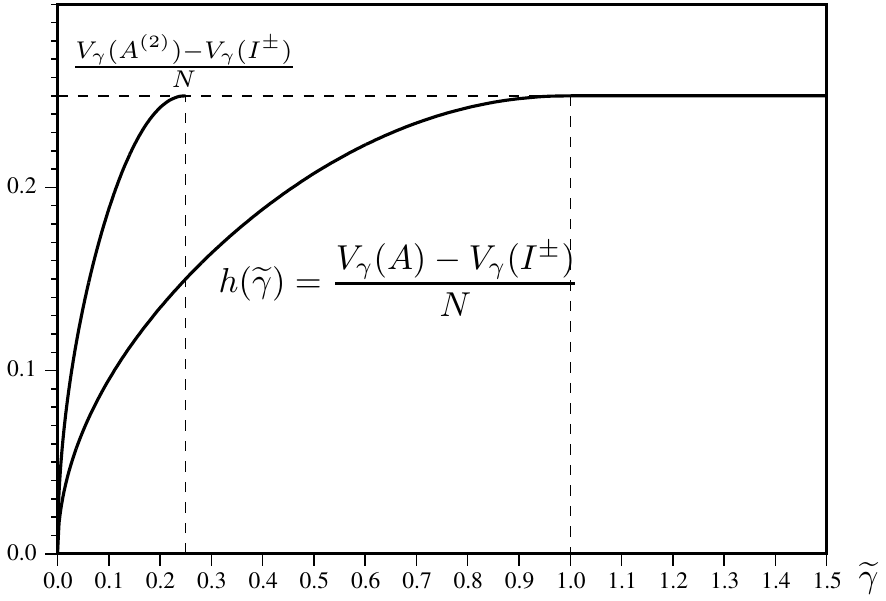}
}
 \figtext{ }
 \caption[]
 {Value of the rescaled potential barrier height
 $h(\gammat)=(V_\gamma(A)-V_\gamma(I^\pm))/N$ as a function of the rescaled
 coupling intensity $\cngt$. For comparison, we also show the
 rescaled barrier
 height $(V_\gamma(A^{(2)})-V_\gamma(I^\pm))/N$ for a stationary point of
 the higher winding number $M=2$.}
\label{fig_largeN2}
\end{figure}

We can now prove a similar result in the desynchronised regime
$\gammat<1$. 
\begin{theorem}
\label{thm_stoch2}
For $\gammat<1$, let 
\begin{equation}
\label{stoch7}
h(\gammat) = \frac{V_\gamma(A(\gammat))-V_\gamma(I^\pm)}N = \frac14 -
\frac1{3(1+\alkappa^2)}
\biggbrak{\frac{2+\alkappa^2}{1+\alkappa^2}  - 2 
\frac{\JE(\alkappa)}{\JK(\alkappa)}} + 
\biggOrder{\frac{\alkappa^2}{N}}\;,
\end{equation}
where $\alkappa=\alkappa(\gammat)$ is defined implicitly
by~\eqref{LargeN3}. Fix an initial condition $x_0\in\cB(I^-,r)$.
Then for any $0<\gammat<1$, and any $\delta>0$, there exists
$N_0(\gammat)$ such that for all $N>N_0(\gammat)$, 
\begin{equation}
\label{stoch8}
\lim_{\sigma\to0}
\bigprobin{x_0}{\e^{(2Nh(\gammat)-\delta)/\sigma^2} < \tau_+ <
\e^{(2Nh(\gammat)+\delta)/\sigma^2}} = 1\;
\end{equation}
and 
\begin{equation}
\label{stoch8b}
\lim_{\sigma\to0} \sigma^2 \log 
\expecin{x_0}{\tau_+} = 2Nh(\gammat)\;.
\end{equation}
Furthermore, let 
\begin{equation}
\label{stoch9}
\tau_A = \tauhit\Bigpar{\bigcup_{g\in G}\cB(gA,r)}\;,
\end{equation}
where $A=A(\gammat)$ satisfies~\eqref{LargeN5} (or~\eqref{largeN7} if $N$
is
odd). Then for any $N>N_0(\gammat)$, 
\begin{equation}
\label{stoch10}
\lim_{\sigma\to0}
\bigpcondin{x_0}{\tau_A < \tau_+}{\tau_+ < \tau_-} = 1\;.
\end{equation}
\end{theorem}

The relations~\eqref{stoch8} and~\eqref{stoch8b} mean that the transition
time between the synchronised states $I^-$ and $I^+$ is of order
$\e^{2Nh(\gammat)/\sigma^2}$, while relation~\eqref{stoch10} implies that
the set of critical configurations is given by the group orbit of $A$ under
$G$.

The large-$N$ limit of the rescaled potential difference $h(\gammat)$ is
shown in~\figref{fig_largeN2}. The limiting function is increasing, with a
discontinuous second-order derivative at $\gammat=1$. For small $\gammat$,
$h(\gammat)$ grows like the square-root of $\gammat$. This is compatible
with the weak-coupling behaviour $h=(1/4+3/2\gamma+\Order{\gamma^2})/N$
obtained in~\cite{BFG06a}, if one takes into account the scaling of
$\gammat$. 

The critical configuration, that is, the configuration with highest energy
reached in the course of the transition from $I^-$ to $I^+$, is any
translate of the configuration $A(\gammat)$ shown in~\figref{fig_symell}a.
If $N$ is even, it has $N/2$ positive and $N/2$ negative coordinates, while
for odd $N$, there are $(N-1)/2$ positive, one vanishing, and $(N-1)/2$
negative coordinates. The sites with positive and negative coordinates are
always adjacent. The potential difference between the $1$-saddles $A$ and
the $2$-saddles $B$ is actually very small, so that transition paths
become less localised as the particle number $N$ increases, reflecting the
fact that the system becomes translation-invariant in the large-$N$ limit.


\newpage
\section{Proofs}
\label{sec_tm}


\subsection{Strategy of the Proof}
\label{sec_tmstrat}


The proof of Theorems~\ref{thm_LargeN1} and~\ref{thm_LargeN2} is based
on the fact that stationary points of the potential satisfy the relation 
\begin{equation}
\label{tm1}
f(x_n) + \frac\cng2 \bigbrak{x_{n+1}-2x_n+x_{n-1}} =
0\;,
\end{equation}
where $f(x)=x-x^3$. As mentioned in Section~\ref{ssec_largeN}, this
relation can be rewritten as a two-dimensional area-preserving twist map
\begin{equation}
\label{tm2}
\begin{split}
x_{n+1} &= x_n + v_{n+1}\;, \\
v_{n+1} &= v_n - 2\gamma^{-1} f(x_n)\;. 
\end{split}
\end{equation}
whose periodic points correspond to stationary points of the potential. In
fact, we are going to analyse a slightly different equivalent map, which
has the advantage to use the symmetries of the model in a more efficient
way. 

The proof is organised as follows:
\begin{itemiz}
\item 	In Section~\ref{sec_tmsym}, we introduce the alternative twist
map, adapted to symmetries.  
\item 	In Section~\ref{sec_tmaa}, we compute the expression of the map in
action--angle variables, taking advantage of the existence of an almost
conserved quantity. 
\item 	In Section~\ref{sec_tmgf}, we compute the generating function of
the map in action--angle variables. This reduces the problem of
finding periodic orbits to a variational problem (which is simpler than
the original one). 
\item 	The main difficulty is that the system is almost degenerate along
the translation mode. In Section~\ref{sec_tmgff}, we introduce Fourier
variables in order to decouple the translation mode from the other, \lq\lq
oscillating\rq\rq\ modes. 
\item 	In Section~\ref{sec_tmgchi}, we deal with the oscillating modes,
by showing with the help of Banach's contraction principle that for each
value of the translation mode, there is exactly one value of the
oscillating modes yielding a stationary point.
\item 	In Section~\ref{sec_tmgstat}, we deal with the translation mode,
by reducing the problem to one dimension, and showing that the generating
function is dominated by its leading Fourier mode in this direction. This
yields the exact number of stationary points. 
\item 	Finally, in Section~\ref{sec_tmgindex} we consider the stability of
the stationary points. 
\end{itemiz}


\subsection{Symmetric Twist Map}
\label{sec_tmsym}

The twist map~\eqref{tm2} does not exploit the symmetries of the
original system in an optimal way. In order to do so, it is more
advantageous to introduce the variable
\begin{equation}
\label{tmsym1}
u_n = \frac{x_{n+1}-x_{n-1}}2
\end{equation}
instead of $v_n$. Then a short computation shows that 
\begin{equation}
\label{tmsym2}
\begin{split}
x_{n+1} &= x_n + u_n - \gamma^{-1} f(x_n)\;,\\
u_{n+1} &= u_n - \gamma^{-1} \bigbrak{f(x_n)+f(x_{n+1})}\;.
\end{split}
\end{equation}
The map $T_1: (x_n,u_n)\mapsto(x_{n+1},u_{n+1})$ is also an
area-preserving twist map. Although it looks more complicated
than the map~\eqref{tm2}, it has the advantage that its inverse is
obtained by changing the sign of $u$, namely 
 \begin{equation}
\label{tmsym3}
\begin{split}
x_n &= x_{n+1} - u_{n+1} - \gamma^{-1} f(x_{n+1})\;,\\
u_n &= u_{n+1} + \gamma^{-1} \bigbrak{f(x_{n+1})+f(x_n)}\;.
\end{split}
\end{equation}
This implies that if we introduce the involutions 
\begin{equation}
\label{tmsym4}
S_1: (x,u) \mapsto (-x,u)
\qquad\text{and}\qquad
S_2: (x,u) \mapsto (x,-u)\;,
\end{equation}
then the map $T_1$ and its inverse are related by 
\begin{equation}
\label{tmsym5}
T_1 \circ S_1 = S_1 \circ (T_1)^{-1}
\qquad\text{and}\qquad
T_1 \circ S_2 = S_2 \circ (T_1)^{-1}\;, 
\end{equation}
as a consequence of $f$ being odd. 
This implies that the images of an orbit of the map under $S_1$ and $S_2$
are also orbits of the map.

For large $N$, it turns out to be useful to introduce the small parameter 
\begin{equation}
\label{tm3}
\eps = \sqrt{\frac2\cng}
= \sqrt{\frac2{{\gamma_1\cngt}}}
= \frac{2\pi}{N\sqrt{\cngt}}
\biggpar{1+\biggOrder{\frac1{N^2}}}\;,
\end{equation}
and the scaled variable $w = u/\eps$. This transforms the map $T_1$
into a map $T_2: (x_n,w_n)\mapsto(x_{n+1},w_{n+1})$ defined by 
\begin{equation}
\label{tm4}
\begin{split}
x_{n+1} &= x_n + \eps w_n - \frac12\eps^2 f(x_n)\;,\\
w_{n+1} &= w_n - \frac12\eps \bigbrak{f(x_n)+f(x_{n+1})}\;.
\end{split}
\end{equation}
$T_2$ is again an area-preserving twist map satisfying 
\begin{equation}
\label{tm4sym}
T_2 \circ S_1 = S_1 \circ (T_2)^{-1}
\qquad\text{and}\qquad
T_2 \circ S_2 = S_2 \circ (T_2)^{-1}\;. 
\end{equation}


\subsection{Action--Angle Variables}
\label{sec_tmaa}

For small $\eps$, we expect the orbits of this map to be close to those of
the differential equation
\begin{equation}
\label{tm5}
\begin{split}
\dot x &= w\;,\\
\dot w &= -f(x)\;,
\end{split}
\end{equation}
which is equivalent to the second-order equation 
$\ddot x = -f(x)$ 
describing the motion of a particle in the\/ \emph{inverted}\/ double-well
potential $-U(x)$, compare~\eqref{heur4}. Solutions
of~\eqref{tm5} can be expressed in terms of Jacobi elliptic
functions. Indeed, the function 
\begin{equation}
\label{tm7}
C(x,w) = \frac12(x^2+w^2) - \frac14 x^4
\end{equation}
being a constant of motion, one sees that $w$ satisfies
\begin{equation}
\label{tm8}
w = \pm \sqrt{(a(C)^2-x^2)(b(C)^2-x^2)/2}\;,
\end{equation} 
where 
\begin{align}
\nonumber
a(C)^2 &= 1 - \sqrt{1-4C}\;,\\
b(C)^2 &= 1 + \sqrt{1-4C}\;.
\label{tm9} 
\end{align}
This can be used to integrate the equation $\dot x=w$, yielding 
\begin{equation}
\label{tm9b}
\frac{b(C)}{\sqrt2}t = \JF \biggpar{\ArcSin\biggpar{\frac
{x(t)}{a(C)}},\alkappa(C)}\;,
\end{equation}
where $\alkappa(C)=a(C)/b(C)$, and $\JF(\phi,\kappa)$ denotes the
incomplete elliptic integral of the first kind. 
The solution of the ODE can be written in terms of standard elliptic
functions as 
\begin{equation}
\label{tm10}
\begin{split}
x(t) &= a(C) \sn\biggpar{\frac{b(C)}{\sqrt2}t,\alkappa(C)}\;,\\
w(t) &= \sqrt{2C} \cn\biggpar{\frac{b(C)}{\sqrt2}t,\alkappa(C)}
\dn\biggpar{\frac{b(C)}{\sqrt2}t,\alkappa(C)}\;.
\end{split}
\end{equation}


We return now to the map $T_2$ defined in~\eqref{tm4}. The explicit
solution
of the continuous-time equation motivates the change of variables $\Phi_1:
(x,w)\mapsto(\varphi,C)$ given by 
\begin{equation}
\label{tm13}
\begin{split}
\varphi &= \frac{\sqrt2}{b(C)} \JF \biggpar{\ArcSin\biggpar{\frac
x{a(C)}},\alkappa(C)}\;,\\
C &= \frac12(x^2+w^2) - \frac14 x^4\;.
\end{split}
\end{equation}
One checks that $\Phi_1$ is again area-preserving. 
The inverse $\Phi_1^{-1}$ is given by 
\begin{equation}
\label{tm11}
\begin{split}
x &= a(C) \sn\biggpar{\frac{b(C)}{\sqrt2}\varphi,\alkappa(C)}\;, \\
w &= \sqrt{2C} \cn\biggpar{\frac{b(C)}{\sqrt2}\varphi,\alkappa(C)}
\dn\biggpar{\frac{b(C)}{\sqrt2}\varphi,\alkappa(C)}\;.
\end{split}
\end{equation}
The elliptic functions $\sn$, $\cn$ and $\dn$ being periodic in their first
argument, with period $4\JK(\alkappa)$, it is convenient to carry out a
further area-preserving change of variables $\Phi_2:
(\varphi,C)\mapsto(\psi,I)$, defined by  
\begin{equation}
\label{tm14}
\psi = \Omega(C) \varphi\;,
\qquad
I = h(C)\;,
\end{equation}
where 
\begin{equation}
\label{tm15}
\Omega(C) = \frac{b(C)}{\sqrt2} \frac{\pi}{2\JK(\alkappa(C))}\;,
\qquad
h(C) = \int_0^C \frac{\6C'}{\Omega(C')}\;.
\end{equation}
Using the facts that $C$ and $b=b(C)$ can be expressed as functions of
$\alkappa=\alkappa(C)$ by $C=\alkappa^2/(1+\alkappa^2)^2$ and 
$b^2=2/(1+\alkappa^2)$, one can check that 
\begin{equation}
\label{tm15a}
h(C) = \frac4{3\pi} 
\frac{(1+\alkappa^2)\JE(\alkappa) - (1-\alkappa^2)\JK(\alkappa)}
{(1+\alkappa^2)^{3/2}}\biggr|_{\alkappa=\alkappa(C)}
\in \biggbrak{0,\frac{2\sqrt2}{3\pi}}\;.
\end{equation}
We denote by $\Phi=\Phi_2\circ\Phi_1$ the transformation
$(x,w)\mapsto(\psi,I)$ and by $T = \Phi\circ T_2\circ\Phi^{-1}$ the
resulting map. 

\begin{prop}
\label{prop_tm1}
The map $T=T(\eps)$ has the form  
\begin{equation}
\label{tm16}
\begin{split}
\psi_{n+1} &= \psi_n + \eps \Omegabar(I_n) + \eps^3 f(\psi_n,I_n,\eps) 
\qquad \pmod{2\pi}\;,\\
I_{n+1} &= I_n + \eps^3 g(\psi_n,I_n,\eps)\;,
\end{split}
\end{equation}
where $\Omegabar(I) = \Omega(h^{-1}(I))$. The functions $f$ and $g$
are $\pi$-periodic in their first argument, and are real-analytic for 
$0\leqs I\leqs h(1/4)-\Order{\eps^3}$. Furthermore, $T$ satisfies the
symmetries 
\begin{equation}
\label{tm16sym}
T \circ \Sigma_1 = \Sigma_1 \circ T^{-1}
\qquad\text{and}\qquad
T \circ \Sigma_2 = \Sigma_2 \circ T^{-1}\;, 
\end{equation}
where $\Sigma_1(\psi,I)=(-\psi,I)$ and $\Sigma_2(\psi,I)=(\pi-\psi,I)$. 
\end{prop}
\begin{proof}
First observe that $\Phi_1$ and $\Phi$ are analytic whenever $(x,w)$ is
such
that $C<1/4$. The map $T$ will thus be analytic whenever $(\psi_n,I_n)$ is
such that $C(x_n,w_n)<1/4$ and $C(x_{n+1},w_{n+1})<1/4$. 
A direct computation shows that 
\begin{equation}
\label{tm16:1}
C(x_{n+1},w_{n+1}) - C(x_n,w_n) = 
\frac{\eps^3}4 \Bigbrak{x_nw_n + 2x_nw_n^3 - 4x_n^3w_n + 3x_n^5w_n} +
\Order{\eps^4}\;.
\end{equation}
This implies that $I_{n+1}=I_n+\Order{\eps^3}$, and allows to determine
$g(\psi,I,0)$. It also shows that $T$ is analytic for
$I_n<h(1/4)-\Order{\eps^3}$. 
Furthermore, writing $a_n=a(C(x_n,w_n))$, we see that~\eqref{tm16:1}
implies
$a_{n+1}-a_n=\Order{\eps^3}$ and similarly for
$b_n$, $\kappa_n$. This yields   
\begin{align}
\nonumber
\varphi(x_{n+1},w_{n+1}) - \varphi(x_n,w_n) 
&= \frac{\sqrt2}{b_n} \int_{x_n/a_n}^{x_{n+1}/a_n}
\frac{\6u}{\sqrt{(1-\kappa_n^2 u^2)(1-u^2)}} + \Order{\eps^3}\\
&= \eps + \Order{\eps^3}\;,
\label{tm16:2}
\end{align}
which implies the expression for $\psi_{n+1}$. We remark that the fact that
$T$ is area-preserving implies the relation  
\begin{equation}
\label{tm16:3}
1 = \dpar{(\psi_{n+1},I_{n+1})}{(\psi_n,I_n)} 
= 1 + \eps^3 \Bigbrak{\sdpar f\psi(\psi,I,0) + \sdpar gI(\psi,I,0)} +
\Order{\eps^4}\;,
\end{equation}
which allows to determine $f(\psi,I,0)$. The fact that $f$ and $g$ are
$\pi$-periodic in their first argument is a consequence of the fact that
$T_2(-x,-w)=-T_2(x,w)$. Finally the relations~\eqref{tm16sym} follow
from the symmetries~\eqref{tm4sym}, with $\Sigma_i=\Phi \circ S_i \circ
\Phi^{-1}$. 
\end{proof}

A perturbation expansion at $I=0$ shows in particular that 
\begin{equation}
\label{tm17a}
\Omegabar(I) = 1-\frac34I+\Order{I^2}\;.
\end{equation}
An important observation is that $\Omegabar(I)$ is a monotonously 
decreasing function, taking values in $[0,1]$. The monotonicity of
$\Omegabar$ makes $T$ a\/ \emph{twist map}\/ for sufficiently small $\eps$,
which has several important consequences on existence of periodic orbits. 

We call\/ \emph{rotation number}\/ of a periodic orbit of period $N$ the
quantity 
\begin{equation}
\label{tm18}
\nu = \frac1{2\pi N} \biggbrak{\sum_{n=1}^N (\psi_{n+1}-\psi_n)
\pmod{2\pi}}\;.
\end{equation}
Note that because of periodicity, $\nu$ is necessarily a rational number of
the form  $\nu=M/N$, for some positive integer $M$.  We denote by
$\T^N_\nu$ the set of points $\psi$ in the torus $\T^N$
satisfying~\eqref{tm18}. It is sometimes more convenient to visualise
$\T^N_\nu$ as the set of real \mbox{$N$-tuples} $(\psi_1,\dots,\psi_N)$
such that 
\begin{equation}
\label{tm19B}
\psi_1 < \psi_2 < \dots < \psi_{N} < \psi_1 + 2\pi N\nu\;.
\end{equation}
In the sequel, we shall use the shorthand\/ \emph{stationary point with
rotation number $\nu$}\/ instead of\/ \emph{stationary point corresponding
to a periodic orbit of rotation number $\nu$}.  

The expression~\eqref{tm16} for $T$ implies that 
\begin{equation}
\label{tm19}
\nu = \frac{\Omegabar(I_0)}{2\pi} \eps + \Order{\eps^2}\;.
\end{equation}
The following properties follow from the Poincar\'e--Birkhoff theorem,
whenever $\eps>0$ is sufficiently small:

\begin{itemiz}
\item	For each positive integer $M$ satisfying 
\begin{equation}
\label{tm20}
M\leqs \frac{N\eps}{2\pi}(1+\Order{\eps})\;,
\end{equation}
the twist map $T$ admits at least two periodic orbits of period $N$ and
rotation number $\nu=M/N$. Note that Condition~\eqref{tm20} is compatible
with the fact that $O$ bifurcates for $\gamma=\gamma_M$,
$M=1,2,\dots,\intpart{N/2}$. 

\item	Any periodic orbit of period $N$ of the map $T$ is of
the form
\begin{equation}
\label{tm21}
\begin{split}
\psi_n &= \psi_0 + 2\pi\nu n + \Order{\eps^2} \;,\\
I_n &= \Omegabar^{-1} \biggpar{\frac{2\pi}\eps \nu} +
\Order{\eps^2}\;,
\end{split}
\end{equation}for some $\psi_0$ and some $\nu=M/N$, where $M$
is a positive integer satisfying~\eqref{tm20}. 
\end{itemiz}

Going back to original variables, we see that these periodic orbits are of
the form 
\begin{equation}
\label{tm22}
\begin{split}
x_n &= a_n \sn\biggpar{\frac{2\JK(\alkappa_n)}\pi \psi_n,\alkappa_n}\;,\\
w_n &= \sqrt{2C_n} \cn\biggpar{\frac{2\JK(\alkappa_n)}\pi
\psi_n,\alkappa_n}
\dn\biggpar{\frac{2\JK(\alkappa_n)}\pi \psi_n,\alkappa_n}\;,\\
\end{split}
\end{equation}
where $a_n=a(C_n)$, $\alkappa_n = \alkappa(C_n)$ and 
\begin{equation}
\label{tm23}
C_n = \Omega^{-1}\biggpar{\frac{2\pi M}{N\eps}} + \Order{\eps}
= \Omega^{-1}\Bigpar{M\sqrt{\cngt}\,} + \Order{\eps}\;.
\end{equation}
This allows in particular to compute the value of the potential at the
corresponding stationary point.

\begin{prop}
\label{prop_tm3}
Let $\eps>0$ be sufficiently small, and let $x^\star$ be a stationary point
of the potential $V_\gamma$, corresponding to an orbit with rotation number
$\nu=M/N$. Then 
\begin{equation}
\label{tm24}
\frac{V_\gamma(x^\star)}{N} 
= -\frac1{3(1+\alkappa^2)} \biggbrak{\frac{2+\alkappa^2}{1+\alkappa^2} - 
2\frac{\JE(\alkappa)}{\JK(\alkappa)}} + \Order{\eps\alkappa^2}\;, 
\end{equation}
where $\alkappa=\alkappa(C)$, and $C$ satisfies $\Omega(C)^2 = M^2\cngt$. 
\end{prop}
\begin{proof}
The expression~\eqref{mod5} for the potential implies that 
\begin{align}
\nonumber
\frac{V_\gamma(x^\star)}N 
&= \frac1N \sum_{n=1}^N \Bigpar{U(x_n)+\frac12 w_n^2 + \Order{\eps^2}} 
= \frac1N \sum_{n=1}^N (w_n^2-C_n + \Order{\eps^2})\\
\nonumber
&= \frac CN \sum_{n=1}^N 
\biggbrak{2 \cn^2\biggpar{\frac{2\JK(\alkappa)}\pi \psi_n,\alkappa}
\dn^2\biggpar{\frac{2\JK(\alkappa)}\pi \psi_n,\alkappa} - 1 
+ \Order{\eps}}\\
&= C \biggbrak{2 \int_0^{2\pi} 
\cn^2\biggpar{\frac{2\JK(\alkappa)}\pi \psi,\alkappa}
\dn^2\biggpar{\frac{2\JK(\alkappa)}\pi \psi,\alkappa)} \6\psi - 1  
+ \Order{\eps}}\;,
\label{tm25}
\end{align}
where $C=\Omega^{-1}(M\sqrt{\cngt})$ and $\alkappa=\alkappa(C)$. The
integral
can then be computed using the change of variables $2\JK(\alkappa)\psi/\pi
=
\JF(\phi,\alkappa)$.  Finally, recall that
$C=\alkappa^2/(1+\alkappa^2)^2$. 
\end{proof}

One can check that $V_\gamma(x^\star)/N$ is a decreasing function of
$\alkappa$, which is
itself a decreasing function of $M^2\cngt$. As a consequence,
$V_\gamma(x^\star)/N$ is
increasing in $M^2\cngt$. This implies in particular that the potential is
larger for larger winding numbers $M$. 

\begin{remark}
\label{rem_tmpot}
The leading term in the expression~\eqref{tm24} for the value of the
potential is the same for all orbits of a given rotation number $\nu$.
Since stationary points of the potential of different index cannot be at
exactly the same height, the difference has to be hidden in the error
terms. In~\cite{BFG06a}, we showed that near the desynchronisation
bifurcation, the potential difference between $1$-saddles and $2$-saddles
is of order $(1-\gammat)^{N/2}$. For large $N$, we expect this difference
to be exponentially small in $1/N$, owing to the fact that near-integrable
maps of a form similar to~\eqref{tm16} are known to admit adiabatic
invariants to that order (cf.~\cite[Theorem~2]{BK0}). 
\end{remark}


\subsection{Generating Function}
\label{sec_tmgf}

In this section, we transform the problem of finding periodic orbits of
the near-integrable map $T$ into a variational problem. 
The fact that $T$ is a twist map allows us to express $I_n$ (and thus
$I_{n+1}$) as a function of $\psi_n$ and $\psi_{n+1}$. A\/
\emph{generating function}\/ of $T$ is a function $G(\psi_n,\psi_{n+1})$
such that 
\begin{equation}
\label{tmgf1}
\sdpar G1(\psi_n,\psi_{n+1}) = -I_n\;, 
\qquad
\sdpar G2(\psi_n,\psi_{n+1}) = I_{n+1}\;. 
\end{equation}
It is known that any area-preserving twist map admits a generating
function,
unique up to an additive constant. Since $T$ depends on the parameter
$\eps$, the generating function $G$ naturally also depends on $\eps$.
However, we will indicate this dependence only when we want to emphasise
it.

\begin{prop}
\label{prop_tmgf}
The map $T$ admits a generating function of the form 
\begin{equation}
\label{tmgf2}
G(\psi_1,\psi_2) = \eps G_0\biggpar{\frac{\psi_2-\psi_1}\eps, \eps} 
+ 2\eps^3 \sum_{p=1}^\infty \Ghat_p\biggpar{\frac{\psi_2-\psi_1}\eps,
\eps} \cos\bigpar{p(\psi_1+\psi_2)}\;,
\end{equation}
where the functions $G_0(u,\eps)$ and $\Ghat_p(u,\eps)$ are real-analytic
for $u>\Order{1/\abs{\log\eps}}$, and satisfy  
\begin{align}
\nonumber
G_0'(u,0) &= \Omegabar^{-1}(u)\;,\\
\Ghat_p\bigpar{u,0}  
&=  \frac1{4p\pi} \int_0^{2\pi} 
g\bigpar{\psi,\Omegabar^{-1}(u),0} \sin(-2p\psi)\,\6\psi\;.
\label{tmgf3}
\end{align}
\end{prop}
\begin{proof}
Fix $(\psi_2,I_2) = T(\psi_1,I_1)$. The fact that
$T(\psi_1+\pi,I_1)=(\psi_2+\pi,I_2)$ implies
\begin{equation}
\label{tmgf3:1}
G(\psi_1+\pi,\psi_2+\pi) = G(\psi_1,\psi_2) + c
\end{equation}
for some constant $c$. If we set
$G(\psi_1,\psi_2)=\Gtilde(\psi_2-\psi_1,\psi_1+\psi_2)$, we thus have 
\begin{equation}
\label{tmgf3:2}
\Gtilde(u,v+2\pi) = \Gtilde(u,v) + c\;.
\end{equation}
This allows us to expand $G$ as a Fourier series 
\begin{equation}
\label{tmgf3:3}
G(\psi_1,\psi_2) 
= \sum_{p=-\infty}^\infty \Gtilde_p(\psi_2-\psi_1,\eps) \e^{\icx
p(\psi_1+\psi_2)} + \frac c{2\pi} (\psi_1+\psi_2)\;.
\end{equation}
Next we note that the symmetry~\eqref{tm16sym} implies 
$T(-\psi_2,I_2) = (-\psi_1,I_1)$, and thus 
\begin{equation}
\label{tmgf3:4}
\sdpar G1(\psi_1,\psi_2) = -\sdpar G2 (-\psi_2,-\psi_1)\;.
\end{equation}
Plugging~\eqref{tmgf3:3} into this relation yields 
\begin{equation}
\label{tmgf3:5}
c = 0
\qquad\text{and}\qquad
\Gtilde_{-p}(u,\eps) = \Gtilde_p(u,\eps)\;,
\end{equation}
which allows to represent $G$ as a real Fourier series as well. 
Computing the derivatives $I_1=-\sdpar G1(\psi_1,\psi_2)$ and
$I_2=\sdpar G2(\psi_1,\psi_2)$ yields 
\begin{equation}
\label{tmgf3:6}
I_2 - I_1 = 2\sum_{p=-\infty}^\infty \icx p \Gtilde_p(\psi_2-\psi_1,\eps) 
\e^{\icx p(\psi_1+\psi_2)}\;,
\end{equation}
which shows in particular that $\Gtilde_p(u,\eps) = \Order{\eps^3}$ for
$p\neq0$,  as a consequence of~\eqref{tm16}.  This implies
$I_1=\Gtilde_0'(\psi_2-\psi_1,\eps) + \Order{\eps^3}$, and thus $u =
\Gtilde_0'(\eps\Omegabar(u)+\Order{\eps^3},\eps)$. Renaming
$\Gtilde_0(u,\eps)=\eps G_0(u/\eps,\eps)$ and $\Gtilde_p(u,\eps)=\eps^3
\Ghat_p(u/\eps,\eps)$ yields~\eqref{tmgf2}. Evaluating~\eqref{tmgf3:6} for
$\eps=0$ and taking the Fourier transform yields the 
expression~\eqref{tmgf3} for $\Ghat_p(u,0)$. 
\end{proof}

The relations~\eqref{tmgf3} allow to determine the expression for the
generating function of the map $T$, given by~\eqref{tm16}. 
In particular, one finds 
\begin{equation}
\label{tmfg4a}
G_0(u,0) = u\Omegabar^{-1}(u) - \Omega^{-1}(u)\;,
\end{equation}
so that 
\begin{align}
\nonumber
G_0(\Omega(C),0) &= h(C)\Omega(C)-C \\
&= -\frac1{3(1+\alkappa^2)} \biggbrak{\frac{2+\alkappa^2}{1+\alkappa^2} - 
2\frac{\JE(\alkappa)}{\JK(\alkappa)}}\;,
\label{tmfg4b}
\end{align}
with $\alkappa=\alkappa(C)$. Note that this quantity is identical with the
leading term in the expression~\eqref{tm24} for the average potential per
site. This indicates that we have chosen the integration constant in the
generating function in such a way that $V_\gamma$ and $G_N$ take the same
value on
corresponding stationary points. 

The main use of the generating function lies in the following fact.
Consider
the $N$-point function 
\begin{equation}
\label{tmgf5}
G_N(\psi_1,\dots,\psi_N) = G(\psi_1,\psi_2) + G(\psi_2,\psi_3) +
\dots + G(\psi_N,\psi_1 + 2\pi N\nu)\;,
\end{equation}
defined on (a subset of) the set $\T^N_\nu$. The defining
property~\eqref{tmgf1} of the generating function implies that for any
periodic orbit of period $N$ of the map $T$, one has 
\begin{equation}
\label{tmgf6}
\dpar{}{\psi_n} G_N(\psi_1,\dots,\psi_N) = -I_n + I_n = 0\;,
\qquad
\text{for $n=1,\dots,N$.}
\end{equation}
In other words, $N$-periodic orbits of $T$ with rotation number $\nu$ are
in
one-to-one correspondence with stationary points of the $N$-point function
$G_N$ on $\T^N_\nu$. 

The symmetries of the original potential imply that the $N$-point
generating function satisfies the following relations on $\T^N_\nu$:
\begin{align}
\nonumber
G_N(\psi_1,\dots,\psi_N) &= G_N(\psi_2,\dots,\psi_N,\psi_1+2\pi N\nu)\;,\\
\nonumber
G_N(\psi_1,\dots,\psi_N) &= G_N(-\psi_N,\dots,-\psi_1)\;,\\
G_N(\psi_1,\dots,\psi_N) &= G_N(\psi_1+\pi,\dots,\psi_N+\pi)\;. 
\label{tmgf7}
\end{align}

At this point, we are in the following situation. We have first transformed
the initial problem of finding the stationary points of the potential
$V_\gamma$
into the problem of finding periodic orbits of the map $T_1$, or,
equivalently, of the map $T$. This problem in turn has been transformed
into the problem of finding the stationary points of $G_N$. Obviously,
the whole procedure is of interest only if the stationary points of $G_N$
are easier to find and analyse than those of $V_\gamma$. This, however, is
the
case here since the $N$-point function is a small perturbation of a
function depending only on the differences $\psi_{n+1}-\psi_n$. In other
words, $G_N$ can be interpreted as the energy of a chain of particles
with a uniform nearest-neighbour interaction, put in a weak external
periodic potential. 


\subsection{Fourier Representation of the Generating Function}
\label{sec_tmgff}

The main difficulty in analysing the stationary points of the $N$-point
generating function $G_N$ comes from the fact that it is almost degenerate
under translations of the form $\psi_n\mapsto\psi_n+c$ $\forall n$. The
purpose of this section is to decouple the translation mode from the other
variables, by introducing Fourier variables. 

We fix $\nu=M/N$.  Any stationary point of $G_N$ on $\T^N_\nu$ admits a
Fourier expansion of the form 
\begin{equation}
\label{tmgff1}
\psi_n = 2\pi\nu n + \sum_{q=0}^{N-1} \psihat_q \omega^{qn}\;,
\end{equation}
where $\omega=\e^{2\pi\icx/N}$, and the Fourier coefficients are
uniquely determined by 
\begin{equation}
\label{tmgff2}
\psihat_q = \frac1N \sum_{n=1}^{N}\omega^{-qn} (\psi_n-2\pi\nu n)
=\cc{\psihat_{-q}}\;.
\end{equation}
Note that $\psihat_q=\psihat_{q+N}$ for all $q$.  Stationary points of
$G_N$ correspond to stationary points of the function $\Gbar_N$, obtained
by expressing $G_N$ in terms of Fourier variables
$(\psihat_0,\dots,\psihat_{N-1})$.  In order to do so, it is convenient
to write  
\begin{align}
\nonumber
\frac{\psi_{n+1}-\psi_n}\eps &= \Delta 
+ \eps^2 \alpha_n(\psihat_1,\dots,\psihat_{N-1}) \;, \\
\psi_n + \psi_{n+1} &= 2\psihat_0 + 2\pi\nu(2n+1) 
+ \eps^2 \beta_n(\psihat_1,\dots,\psihat_{N-1}) \;, 
\label{tmgff3}
\end{align}
where $\Delta=2\pi\nu/\eps$ and 
\begin{align}
\nonumber
\alpha_n(\psihat_1,\dots,\psihat_{N-1})
&= \frac1{\eps^3} \sum_{q=1}^{N-1} \psihat_q (\omega^q-1) \omega^{qn}\;,\\
\beta_n(\psihat_1,\dots,\psihat_{N-1})
&= \frac1{\eps^2} \sum_{q=1}^{N-1} \psihat_q (\omega^q+1) \omega^{qn}\;. 
\label{tmgff4}
\end{align}
Note that $\alpha_n$ is of order $1$ in $\eps$ for any stationary point
because of the expression~\eqref{tm16} of the twist map. Taking the inverse
Fourier transform shows that $\abs{\psihat_q(\omega^q-1)}=\Order{\eps^3}$
and $\abs{\psihat_q}=\Order{\eps^2}$ for $q\neq0$, and thus $\beta_n$ is
also of order $1$. 

\begin{table}[b]
\begin{center}
\begin{tabular}{|l|l|l|l|}
\hline
\vrule height 14pt depth 6pt width 0pt
$R$ & $x_j \mapsto x_{j+1}$ & $\psi_n \mapsto \psi_{n+1}$ 
& $\psihat_q \mapsto \omega^q \psihat_q + 2\pi\nu\delta_{q0}$ \\
\vrule height 8pt depth 6pt width 0pt
$CS$ & $x_j \mapsto -x_{N+1-j}$ & $\psi_n \mapsto -\psi_{N+1-n}$ 
& $\psihat_q \mapsto -\omega^{-q}\psihat_{N-q} - 2\pi\nu (N+1)\delta_{q0}$
\\
\vrule height 6pt depth 8pt width 0pt
$C$ & $x_j \mapsto -x_j$ & $\psi_n \mapsto \psi_n+\pi$ 
& $\psihat_q \mapsto \psihat_q + \pi\delta_{q0}$ \\
\hline 
\end{tabular} 
\end{center}
\caption[]
{Effect of some symmetries on original variables, angle variables,  and
Fourier variables.}
\label{table_Fourier_gen}
\end{table}

Expressing $G_N$ in Fourier variables yields the function 
\begin{equation}
\label{tmgff5}
\Gbar_N(\psihat_0,\dots,\psihat_{N-1}) 
= \sum_{p=-\infty}^\infty \e^{2\icx p\psihat_0}
g_p(\psihat_1,\dots,\psihat_{N-1})\;,
\end{equation}
where (we drop the $\eps$-dependence of $G_0$ and $\Ghat_p$) 
\begin{align}
\nonumber
g_0(\psihat_1,\dots,\psihat_{N-1})
&= \eps \sum_{n=1}^N G_0(\Delta+\eps^2 \alpha_n)\;,\\
g_p(\psihat_1,\dots,\psihat_{N-1})
&= \eps^3 \sum_{n=1}^N \Ghat_p(\Delta+\eps^2 \alpha_n) \omega^{pM(2n+1)}
\e^{\icx \eps^2 p \beta_n}
\qquad \text{for $p\neq0$\;.}
\label{tmgff6}
\end{align}

We now examine the symmetry properties of the Fourier coefficients $g_p$. 
Table~\ref{table_Fourier_gen} shows how the Fourier variables transform
under some symmetry transformations, where we only consider transformations
leaving $\T^N_\nu$ invariant. As a consequence, the first two symmetries
in~\eqref{tmgf7} translate into 
\begin{align}
\nonumber
g_p(\psihat_1,\dots,\psihat_{N-1}) 
&= \omega^{2pM} g_p(\omega\psihat_1,\dots,\omega^{N-1}\psihat_{N-1})\;,\\
g_p(\psihat_1,\dots,\psihat_{N-1}) 
&= \omega^{-2pM}
g_{-p}(-\omega^{N-1}\psihat_{N-1},\dots,-\omega\psihat_1)\;. 
\label{tmgff7}
\end{align}
We now introduce new variables $\chi_q$, $q\neq0$, defined by 
\begin{equation}
\label{tmgff9}
\chi_q = -\icx \omega^{-q\psihat_0/2\pi\nu} \psihat_q
= -\cc{\chi_{-q}}\;.
\end{equation} 
The $\chi_q$ are defined in such a way that they are real for stationary
points satisfying, in original variables, the symmetry $x_j=-x_{n_0-j}$ for
some $n_0$. For later convenience, we prefer to consider $q$ as belonging
to 
\begin{equation}
\label{tmgff9B}
\cQ=\biggset{-\biggintpart{\frac{N-1}2},\dots,\biggintpart{\frac N2}}
\setminus\bigset{0}
\end{equation}
rather than $\set{1,\dots,N-1}$. We set $\chi=\set{\chi_q}_{q\in\cQ}$ and 
\begin{align}
\nonumber
\Gtilde_N(\psihat_0,\chi) 
&=
\Gbar_N(\psihat_0,\set{\psihat_q=\icx\omega^{q\psihat_0/2\pi\nu}\chi_q}_{
q\in\cQ}) \\
&= \sum_{p=-\infty}^\infty \e^{2\icx p\psihat_0}
\tilde g_p(\psihat_0,\chi)\;,
\label{tmgff10}
\end{align}
where
\begin{equation}
\label{tmgff11}
\tilde g_p(\psihat_0,\chi) = 
g_p(\set{\psihat_q=\icx\omega^{q\psihat_0/2\pi\nu}\chi_q}_{q\in\cQ})\;.
\end{equation}

\begin{lemma}
\label{lem_tmgff1}
The function $\Gtilde_N(\psihat_0,\chi)$ is
$2\pi\nu$-periodic in its first argument. 
\end{lemma}
\begin{proof}
By~\eqref{tmgff7}, we have 
\begin{equation}
\label{tmgff11:1}
\tilde g_p(\psihat_0+2\pi\nu,\chi)
= \omega^{-2pM} \tilde g_p(\psihat_0,\chi)
\end{equation}
Since $\e^{2\icx p\cdot 2\pi\nu}=\omega^{2pM}$, replacing $\psihat_0$ by
$\psihat_0+2\pi\nu$ in~\eqref{tmgff10} leaves $\Gtilde_N$ invariant. 
\end{proof}

Since $\Gtilde_N$ also has period $\pi$, it has in fact period 
\begin{equation}
\label{tmgff12}
\frac{\pi}N K\;,
\qquad
K = \gcd(N,2M)\;.
\end{equation}

Our strategy now proceeds as follows:
\begin{enum}
\item	Show that for each $\psihat_0$, and sufficiently small $\eps$, the
equations $\tdpar{\Gtilde_N}{\chi_q}=0$, $q\in\cQ$, admit exactly one
solution $\chi=\chi^\star(\smash{\psihat_0})$. This is done in 
Section~\ref{sec_tmgchi} with the help of Banach's fixed-point theorem. 
\item	Show that for $\chi=\chi^\star(\psihat_0)$, the equation 
$\tdpar{\Gtilde_N}{\psihat_0}=0$ is satisfied by exactly $4N/K$ values of 
$\smash{\psihat_0}$. This is done in Section~\ref{sec_tmgstat} by
estimating the Fourier coefficients of $\tdpar{\Gtilde_N}{\psihat_0}$ with
the help of complex analysis. 
\end{enum}


\subsection{Uniqueness of $\mathbf{\chi}$}
\label{sec_tmgchi}

In this section, we show that the equations 
\begin{equation}
\label{tmgchi001}
\dpar{\Gtilde_N}{\chi_q}=0\;,
\qquad
q\in\cQ\;,
\end{equation} 
admit exactly one solution $\chi=\chi^\star(\smash{\psihat_0})$ for each
value of $\psihat_0$. The proof is based on a standard fixed-point
argument: First we show in Lemma~\ref{lem_tmgff4} that~\eqref{tmgchi001} is
equivalent to the fixed-point equation $\rho=\cT\rho$ for a quantity $\rho$
related to $\chi$. Then we show in Proposition~\ref{cor_tmgff2} that $\cT$
is contracting in an appropriate norm, provided $\eps$ is sufficiently
small. 

It is useful to introduce the scaled variables 
\begin{equation}
\label{tmgffB1}
\rho_q = \rho_q(\chi) = -\frac2{\eps^3} \chi_q \sin(\pi q/N)
\end{equation}
and the function $\Gamma^{(a,b)}_{\ell}(\rho)$,
$\rho=\set{\rho_q}_{q\in\cQ}$, defined for $\ell\in\Z$ and $a,
b\geqs 0$ by 
\begin{equation}
\label{tmgffB2}
\Gamma^{(a,b)}_{\ell}(\rho) = 
\sum_{\substack{q_1,\dots,q_a\in\cQ \\ q'_1,\dots,q'_b\in\cQ}} 
\indexfct{\sum_i q_i + \sum_j q'_j = \ell}
\prod_{i=1}^a \rho_{q_i}
\prod_{j=1}^b \frac{\eps}{\tan(\pi q'_j/N)} \rho_{q'_j}\;.
\end{equation}
By convention, any term in the sum for which $q'_j=N/2$ for some $j$ is
zero, that is, we set $1/\tan(\pi/2)=0$.  A few elementary properties
following immediately from this definition are:
\begin{itemiz}
\item	$\Gamma^{(0,0)}_{\ell}(\rho) = \delta_{\ell0}$;
\item	$\Gamma^{(a,b)}_{\ell}(\rho) = 0$ for $\abs{\ell} > (a+b)N/2$;
\item	If $\rho_q=0$ for $q\not\in K\Z$, then
$\Gamma^{(a,b)}_{\ell}(\rho)=0$
for $\ell\not\in K\Z$; 
\item	If $\rho'_q=\rho_{-q}$ for all $q$, then
$\Gamma^{(a,b)}_{\ell}(\rho')=(-1)^b\Gamma^{(a,b)}_{-\ell}(\rho)$. 
\end{itemiz}

The following result states that the conditions~\eqref{tmgchi001} are
equivalent to a fixed-point equation. 

\begin{lemma}
\label{lem_tmgff4}
Let 
\begin{equation}
\label{tmgffA3}
H_{p,q}(\Delta) = \Ghat_p'(\Delta) 
- \frac{\eps p \Ghat_p(\Delta)}{\tan(\pi q/N)}\;,
\end{equation}
with the convention that $H_{p,N/2}(\Delta) = \Ghat_p'(\Delta)$. 
Then the stationarity conditions~\eqref{tmgchi001} are fulfilled if and
only if $\rho=\rho(\chi)$ satisfies the fixed-point equation 
\begin{equation}
\label{tmgffB3}
\rho = \cT \rho = \rho^{(0)} + \Phi(\rho,\eps)\;,
\end{equation}
where the leading term is given by 
\begin{equation}
\label{tmgffB4}
\rho^{(0)}_q = 
\begin{cases}
\displaystyle
\frac{1}{G_0''(\Delta)} \sum_{k\in\Z \colon kN+q\in2M\Z}
(-1)^{k+1} \e^{\icx k\psihat_0 N/M} 
H_{(kN+q)/2M,q}(\Delta)
&\text{if $q\in K\Z$\;,}\\
0
&\text{if $q\not\in K\Z$\;,}
\end{cases}
\end{equation}
and the remainder is given by 
$\Phi_q(\rho,\eps)=\Phi^{(1)}_q(\rho,\eps)+\Phi^{(2)}_q(\rho,\eps)$, with 
\begin{align}
\nonumber
\Phi^{(1)}_q(\rho,\eps)
&= \frac{1}{G_0''(\Delta)} \sum_{k\in\Z} (-1)^{k+1} \e^{\icx k\psihat_0
N/M} 
\sum_{a\geqs 1} \frac{\eps^{2a}}{(a+1)!}
G_0^{(a+2)}(\Delta) 
\Gamma^{(a+1,0)}_{kN+q}(\rho)\;, \\
\Phi^{(2)}_q(\rho,\eps)
&= \frac{1}{G_0''(\Delta)} \sum_{k\in\Z} (-1)^{k+1} \e^{\icx k\psihat_0
N/M} 
 \sum_{a+b\geqs1} \frac{\eps^{2(a+b)}}{a!b!}
\sum_{p\neq0} H^{(a)}_{p,q}(\Delta) p^b
\Gamma^{(a,b)}_{kN-2pM+q}(\rho)\;.
\label{tmgffB5}
\end{align}
\end{lemma}

The proof is a straightforward but lengthy computation, which we postpone
to Appendix~\ref{sec_prtech}. 

Note the following symmetries, which follow directly from the definition of
$\rho^{(0)}$ and the properties of $\smash{\Gamma^{(a,b)}_\ell}$:
\begin{itemiz}
\item	For all $q\in\cQ$, $\rho^{(0)}_{-q}=\rho^{(0)}_q$, because
$H_{-p,-q}(\Delta)=H_{p,q}(\Delta)$, and thus $\rho^{(0)}_q\in\R$;
\item	If $\rho_q=0$ for $q\not\in K\Z$, then $\Phi_q(\rho,\eps)=0$
for $q\not\in K\Z$; 
\item	If $\rho'_q=\rho_{-q}$ for all $q$, then
$\Phi_q(\rho',\eps)=\Phi_{-q}(\rho,\eps)$;
\end{itemiz}

\begin{remark}
\label{rem_tmgffB}
The condition $kN+q\in 2M\Z$, appearing in the definition of $\rho^{(0)}$, 
can only be fulfilled if $q\in N\Z+2M\Z = K\Z$. If this is the case, set
$N=nK$, $2M=mK$, $q=\ell K$, with $n$ and $m$ coprime. Then the condition 
becomes $mp-kn=\ell$.  By Bezout's theorem, the general solution is given
in terms of any particular solution $(p_0,k_0)$ by 
\begin{equation}
\label{tmgff17:3}
p = p_0 + n t\;, 
\qquad
k = k_0 + m t\;,
\qquad
t \in \Z\;.
\end{equation}
Thus there will be exactly one $p$ with $2\abs{p}< n$. 
If $N$ is very large, and $M$ is fixed, then $n=N/K$ is also very large.
Since the $\Ghat_p(\Delta)$, being Fourier coefficients of an analytic
function, decrease exponentially fast in $\abs{p}$, the sum
in~\eqref{tmgffB4} will be dominated by the term with the lowest possible
$\abs{p}$. 
\end{remark}

We now introduce the following weighted norm on $\C^\cQ$:
\begin{equation}
\label{tmgff20}
\norm{\rho}_\lambda = \sup_{q\in\cQ}
\e^{\lambda\abs{q}/2M}\abs{\rho_q}\;,
\end{equation}
where $\lambda>0$ is a free parameter. One checks that the functions
$G_0(\Delta)$ and $\Ghat_p(\Delta)$ are analytic for
$\re\Delta>\Order{1/\log\abs{\eps}}$. Thus it follows from Cauchy's theorem
that there exist positive constants $L_0$,
$r<\Delta-\Order{1/\log\abs{\eps}}$ 
and $\lambda_0$ such that 
\begin{equation}
\label{tmgff21}
\abs{G_0^{(a)}(\Delta)} \leqs L_0 \frac{a!}{r^a}
\qquad\text{and}\qquad
\abs{\Ghat_p^{(a)}(\Delta)} \leqs L_0 \frac{a!}{r^a} \e^{-\lambda_0\abs{p}}
\end{equation}
for all $a\geqs0$ and $p\in\Z$. For sufficiently small $\eps$, it is
possible to choose $r=\Delta/2$. 

\begin{prop}
\label{cor_tmgff2}
There exists a numerical constant $c_1>0$ such that for any $\Delta>0$, any
$\lambda<\lambda_0$ and any $R_0 >
c_1L_0\brak{\Delta\abs{G_0''(\Delta)}}^{-1}$, there is a strictly positive 
$\eps_0=\eps_0(\Delta,\lambda,\lambda_0,R_0)$ such that for all
$\eps<\eps_0$, the map $\cT$ admits a unique fixed point $\rho^\star$ in
the ball
$\cB_\lambda(0,R_0)=\setsuch{\rho\in\C^\cQ}{\norm{\rho}_\lambda<R_0}$.
Furthermore, the fixed point satisfies 
\begin{itemiz}
\item	$\rho^\star_q = 0$ whenever $q\not\in K\Z$;
\item	$\rho^\star_{-q}=\rho^\star_q$, and thus $\rho^\star_q\in\R$ for
all
$q$.
\end{itemiz}
\end{prop}

The proof is again a straightforward but lengthy computation, so we
postpone it to Appendix~\ref{sec_prtech}.

A direct consequence of this result is that for any $\psihat_0$, and
sufficiently small $\eps$, there is a unique
$\rho^\star=\rho^\star(\psihat_0)$ (and thus a unique
$\chi^\star(\psihat_0)$) satisfying the equations
$\tdpar{\Gtilde_N}{\chi_q}=0$ for all $q\in\cQ$. Indeed, we take $R_0$
sufficiently large that our a priori estimates on the $\chi_q$ imply that
$\rho\in\cB_0(0,R_0)$. Then it follows that $\rho$ is unique. Furthermore,
for any $\lambda<\lambda_0$, making $\eps$ sufficiently small we obtain an
estimate on $\norm{\rho^\star}_\lambda$. 


\subsection{Stationary Values of $\mathbf{\psihat_0}$}
\label{sec_tmgstat}

We now consider the condition $\tdpar{\Gtilde_N}{\psihat_0}=0$. As pointed
out at the end of Section~\ref{sec_tmgff}, $\Gtilde_N(\psihat_0,\chi)$  is
a $\pi K/N$-periodic function of $\psihat_0$.  For the same reasons,
$\chi^\star(\psihat_0)$ is also $\pi K/N$-periodic. Hence it follows that
the function $\psihat_0\mapsto \Gtilde_N(\psihat_0,\chi^\star(\psihat_0))$
has the same period as well, and thus admits a Fourier series of the form 
\begin{equation}
\label{tmgstat1}
\Gtilde_N(\psihat_0,\chi^\star(\psihat_0))  
= \sum_{k=-\infty}^\infty \hat g_k \e^{2\icx k \psihat_0 N/K}\;, 
\end{equation}
with Fourier coefficients 
\begin{equation}
\label{tmgstat2}
\hat g_k = \frac 1{2\pi} \int_0^{2\pi} \sum_{p=-\infty}^\infty 
\e^{2\icx(p-kN/K)\psihat_0} \tilde
g_p(\psihat_0,\chi^\star(\psihat_0))\,\6\psihat_0
\end{equation}
(we have chosen $[0,2\pi]$ as interval of integration for later
convenience). 
Using the change of variables $\psihat_0\mapsto-\psihat_0$ in the integral,
and the various symmetries of the coefficients (in
particular~\eqref{tmgff7}), one checks that  $\hat g_{-k}=\hat g_{k}$.
Therefore~\eqref{tmgstat1} can be rewritten in real form as  
\begin{equation}
\label{tmgstat3}
\Gtilde_N(\psihat_0,\chi^\star(\psihat_0))  
= \hat g_0 + 2\sum_{k=1}^\infty \hat g_k \cos\bigpar{2 k \psihat_0 N/K}\;. 
\end{equation}
Now $\tdpar{\Gtilde_N}{\psihat_0}$ vanishes if and only if the total
derivative of $\Gtilde_N(\psihat_0,\chi^\star(\psihat_0))$ with respect to
$\psihat_0$ is equal to zero. This function obviously vanishes for
$\psihat_0=\ell\pi K/2N$, $\ell=1,\dots,4N/K$, and we have to show that
these are the only roots. 

We first observe that the Fourier coefficients $\hat g_k$ can be
expressed directly in terms of the generating function~\eqref{tmgf2},
written in the form 
\begin{equation}
\label{tmgstat4}
\Gtilde(u,v,\eps) = \eps G_0(u,\eps) 
+ \eps^3 \sum_{p\neq0} \Ghat_p(u,\eps) \e^{2\icx pv}\;. 
\end{equation}
In the sequel, $\alpha^\star_n(\psihat_0)$ and $\beta^\star_n(\psihat_0)$
denote the quantities introduced in~\eqref{tmgff4}, evaluated at
$\psihat_q=\icx\omega^{q\psihat_0/2\pi\nu}\chi^\star_q(\psihat_0)$. 

\begin{lemma}
\label{lem_tmgstat1}
The Fourier coefficients $\hat g_k$ are given in terms of the
generating function by 
\begin{equation}
\label{tmgstat5}
\hat g_k = \frac N{2\pi} \int_0^{2\pi}
\e^{-2\icx k\psihat_0N/K} 
\Lambda_0(\psihat_0) \,\6\psihat_0\;,
\end{equation}
where
\begin{equation}
\label{tmgstat5B}
\Lambda_0(\psihat_0) = 
\Gtilde \Bigpar{\Delta+\eps^2\alpha^\star_0(\psihat_0), 
\psihat_0+\pi\nu+\frac12\eps^2\beta^\star_0(\psihat_0),\eps}\;.
\end{equation}
\end{lemma}
\begin{proof}
The coefficient $\hat g_k$ can be rewritten as 
\begin{equation}
\label{tmgstat5:1}
\hat g_k = \frac1{2\pi} \sum_{n=1}^N \int_0^{2\pi}
\e^{-2\icx k\psihat_0N/K} \Lambda_n(\psihat_0) \,\6\psihat_0\;,
\end{equation}
where 
\begin{align}
\nonumber
\Lambda_n(\psihat_0) ={}& \eps G_0(\Delta+\eps^2\alpha^\star_n(\psihat_0))
\\
&{}+ \eps^3\sum_{p\neq0} p \e^{2\icx p\psihat_0} 
\Ghat_p(\Delta+\eps^2\alpha^\star_n(\psihat_0)) \omega^{pM(2n+1)}
\e^{\icx \eps^2 p \beta^\star_n(\psihat_0)} \;.
\label{tmgstat5:2}
\end{align}
Using the periodicity of $\chi^\star$, one finds that 
$\alpha^\star_n(\psihat_0+2\pi\nu) = \alpha^\star_{n+1}(\psihat_0)$ and 
similarly for $\beta^\star_n$, which implies  
$\Lambda_n(\smash{\psihat_0}) = \Lambda_0(\smash{\psihat_0}+2\pi\nu n)$.
Inserting this into~\eqref{tmgstat5:1} and using the change of variables
$\psihat_0\mapsto\psihat_0-2\pi\nu$ in the $n$th summand allows to express 
$\hat g_k$ as the $(2kN/K)$th Fourier coefficient of $\Lambda_0$.  Finally,
$\Lambda_0(\psihat_0)$ can also be written in the form~\eqref{tmgstat5B}. 
\end{proof}

Relation~\eqref{tmgstat5} implies that the $\hat g_k$ decrease
exponentially
fast with $k$, like $\e^{-2\lambda_0kN/K}$. Hence the Fourier 
series~\eqref{tmgstat3} is dominated by the first two terms, provided $N$
is large enough. In order to obtain the existence of exactly $4N/K$
stationary points, it is thus sufficient to prove that $\hat g_1$ is
also bounded below by a quantity of order $\e^{-2\lambda_0N/K}$. 

\begin{prop}
\label{prop_tmgstat2}
For any $\Delta>0$, there exists $\eps_1(\Delta)>0$ such that
whenever $\eps<\eps_1(\Delta)$, 
\begin{equation}
\label{tmgstat6a}
\sign(\hat g_1) = (-1)^{1+2M/K}\;.
\end{equation}
Furthermore, 
\begin{equation}
\label{tmgstat6b}
\frac{\abs{\hat g_k}}{\abs{\hat g_1}} 
\leqs \exp\biggset{-\frac{3k-5}4 \lambda_0(\Delta)\frac NK}
\qquad
\forall k \geqs 2\;,
\end{equation}
where $\lambda_0(\Delta)$ is a monotonously increasing function of
$\Delta$, satisfying $\lambda_0(\Delta)=\sqrt2\pi\Delta+\Order{\Delta^2}$
as $\Delta\searrow0$, and diverging logarithmically as $\Delta\nearrow1$.
\end{prop}
\begin{proof}
First recall that $\Delta=2\pi\nu/\eps=2\pi M/N\eps$, where $M$ is fixed.
Thus taking $\eps$ small for given $\Delta$ automatically yields a large
$N$. 
Combining the expression~\eqref{tm16} for the twist map and the defining
property~\eqref{tmgf1} of the generating function with the relations
$u=(\psi_{n+1}-\psi_n)/\eps$ and $v=\psi_n+\psi_{n+1}$, one obtains the
relation 
\begin{equation}
\label{tmgstat6:1}
\sdpar{\Gtilde}v(u,v,\eps) = \frac{\eps^3}2
\Bigbrak{g\Bigpar{\tfrac12(v-\eps
u), \Omegabar^{-1}(u),\eps} + \Order{\eps^2}}\;.
\end{equation}
It follows from~\eqref{tm16:1} and the
definition~\eqref{tm15} of $h(C)$ that 
\begin{align}
\nonumber
g(\psi,I,0) &= \frac1{\Omegabar(I)} \frac{xw}4 \brak{1+2w^2-4x^2+3x^4}\\
&= \frac1{\Omegabar(I)} \frac{xw}4 \brak{1+4C-6x^2+4x^4}\;,
\label{tmgstat6:2}
\end{align}
where $x$ and $w$ have to be expressed as functions of $\psi$ and $I$
via~\eqref{tm13} and~\eqref{tm14}. 
In particular, we note that 
\begin{equation}
\label{tmgstat6:3}
w = \frac{\sqrt{2C}}a \frac{\pi}{2\JK(\kappa)} \dtot x\psi 
= {\Omegabar(I)} \dtot x\psi\;, 
\end{equation}
where we used~\eqref{tm15} again. This allows us to write
\begin{equation}
\label{tmgstat6:4}
g(\psi,I,0) = \frac18 \dtot{}{\psi}
\biggbrak{(1+4C)x^2 - 3 x^4 + \frac43 x^6}\;. 
\end{equation}
A similar argument would also allow to express the first-order term in
$\eps$ of $g(\psi,I,\eps)$ as a function of $x=x(\psi,I)$. Also note the
equality   
\begin{equation}
\label{tmgstat6:5}
\Delta
= \Omegabar(I) + \Order{\eps}
= \frac{\pi b(C)}{2\sqrt2\JK(\kappa)} + \Order{\eps}
= \frac{1}{\sqrt{1+\kappa^2}\JK(\kappa)} + \Order{\eps}\;,
\end{equation}
which follows from the relation~\eqref{tm19} between $\nu$ and
$\Omegabar(I)$.  

The properties of elliptic functions imply that for fixed $I$, $\psi\mapsto
x(\psi,I)$ is periodic in the imaginary direction, with period $2\lambda_0
= \pi\JK(\sqrt{1-\kappa^2})/\JK(\kappa)$, and has poles located in
$\psi=n\pi + (2m+1)\icx\lambda_0$, $n,m\in\Z$. As a consequence,  the
definition of the map $T = \Phi\circ T_2\circ\Phi^{-1}$ implies in
particular that $g(\psi,I,\eps)$ is a meromorphic function of $\psi$, with
poles at the same location, and  satisfying
$g(\psi+2\icx\lambda_0,I,\eps)=g(\psi,I,\eps)$.  These properties yield
informations on periodicity and location of poles for
$\Lambda_0(\psihat_0)$, in particular
$\Lambda_0(\psihat_0+2\icx\lambda_0)=\Lambda_0(\psihat_0)+\Order{\eps^2}$. 

\begin{figure}
\centerline{\includegraphics*[clip=true,height=40mm]{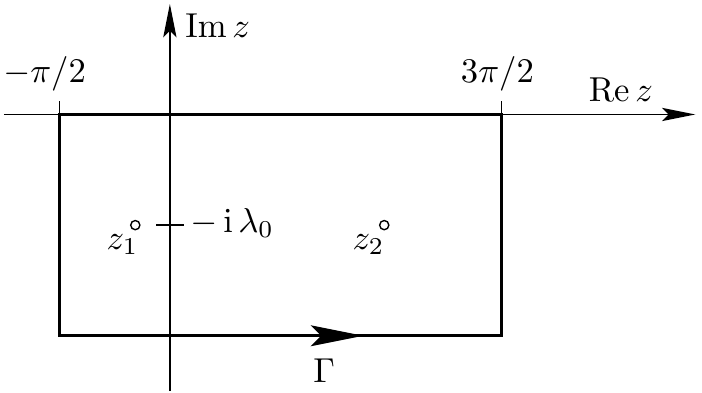}}
 \figtext{ }
 \caption[]
 {The integration contour $\Gamma$ used in the
integral~\eqref{tmgstat6:6}.}
\label{fig_contour}
\end{figure}

Let $\Gamma$ be a rectangular contour with vertices in $-\pi/2$,
$-\pi/2-2\icx\lambda_0$, $3\pi/2-2\icx\lambda_0$, and $3\pi/2$, followed in
the anticlockwise direction (\figref{fig_contour}), and consider the
contour integral
\begin{equation}
\label{tmgstat6:6}
J = \frac1{2\pi} \oint_\Gamma \e^{-2\icx k z N/K} \Lambda_0(z)\,\6z\;.
\end{equation}
The contributions of the integrals along the vertical sides of the
rectangle
cancel by periodicity. Therefore, by Lemma~\ref{lem_tmgstat1} and the
approximate periodicity of $\Lambda_0$ in the imaginary direction, 
\begin{equation}
\label{tmgstat6:7}
J = -\frac1N \Bigbrak{\hat g_k 
- \e^{-2k\lambda_0N/K}\Bigpar{\hat g_k+ \Order{N\eps^5}}}\;.
\end{equation}  
On the other hand, the residue theorem yields 
\begin{equation}
\label{tmgstat6:8}
J = 2\pi\icx \sum_{z_j} 
\e^{-2\icx kz_jN/K}\res(\Lambda_0(z),z_j)\;,
\end{equation}  
where the $z_j$ denote the poles of the function $\Lambda_0(z)$, lying
inside $\Gamma$. There are two such poles, located in
$z_1=-\icx\lambda_0-\pi\nu+\eps\Delta+\Order{\eps^2}$, and $z_2=z_1+\pi$,
and they both yield the same contribution, of order
$\e^{-\lambda_0kN(1+\Order{\eps^2})/K}$, to the sum.
Comparing~\eqref{tmgstat6:7} and~\eqref{tmgstat6:8} shows that $\hat g_k/N$
is of the same order. Finally, the leading term of $\hat g_1$ can be
determined explicitly using~\eqref{tmgstat6:4} and Jacobi's
expression~\eqref{ell9} for the Fourier coefficients of powers of elliptic
functions, and is found to have sign $(-1)^{1+2M/K}$ for sufficiently large
$N$. Choosing $\eps$ small enough (for fixed $\eps N$) guarantees that
$\hat g_1$ dominates all $\hat g_k$ for $k\geqs2$. 
\end{proof}

\begin{cor}
\label{cor_tmgstat}
For $\eps<\eps_1(\Delta)$, the $N$-point generating function $\Gtilde_N$
admits exactly $4N/K$ stationary points, given by $\psihat_0=\ell\pi K/2N$,
$\ell=1,\dots,4N/K$, and $\chi=\chi^\star(\psihat_0)$. 
\end{cor}
\begin{proof}
In the points $\psihat_0=\ell\pi K/2N$,  the derivative of the function
$\psihat_0\mapsto\Gtilde_N(\psihat_0,\chi^\star(\psihat_0))$ vanishes,
while its
second derivative is bounded away from zero, as a consequence
of Estimate~\eqref{tmgstat6b}. Thus these points are simple roots of the
first 
derivative, which is bounded away from zero everywhere else.  
\end{proof}


\subsection{Index of the Stationary Points}
\label{sec_tmgindex}

We finally examine the stability type of the various stationary points, by
first determining their index as stationary points of the $N$-point
generating function $\Gbar_N$, and then examining how this translates into
their index as stationary points of the potential $V_\gamma$. 

\begin{prop}
\label{prop_tmgindex}
Let $(\psihat_0,\chi^\star(\psihat_0))$ be a stationary point of
$\Gtilde_N$ with rotation number $\nu=N/M$. Let
$x^\star=x^\star(\psihat_0)$ be the corresponding stationary point of the
potential $V_\gamma$, and let $K=\gcd(N,2M)$.  
\begin{itemiz}
\item	If $2M/K$ is odd, then the points $x^\star(0)$, $x^\star(K\pi/N)$,
\dots are saddles of even index of $V_\gamma$, while the points
$x^\star(K\pi/2N)$,
$x^\star(3K\pi/2N)$, \dots are saddles of odd index of $V_\gamma$.
\item	If $2M/K$ is even, then the points $x^\star(0)$, $x^\star(K\pi/N)$,
\dots are saddles of odd index of $V_\gamma$, while the points
$x^\star(K\pi/2N)$,
$x^\star(3K\pi/2N)$, \dots are saddles of even index of $V_\gamma$.
\end{itemiz}
\end{prop} 
\begin{proof}
We first determine the index of $(\psihat_0,\chi^\star(\psihat_0))$ as
stationary point of $\Gtilde_N$. Using the fact that $G_0''(\Delta)$ is
negative ($\smash{\Omegabar^{-1}}(\Delta)$ being decreasing), one sees that
the Hessian matrix of $\smash{\Gtilde_N}$ is a small perturbation of a
diagonal
matrix with $N-1$ negative eigenvalues. The $N$th eigenvalue, which
corresponds to translations of $\psihat_0$, has the same sign as the second
derivative of $\psihat_0\mapsto\Gtilde_N(\psihat_0,\chi^\star(\psihat_0))$,
which is equal to $(-1)^{2M/K}\sign\cos(2\psihat_0 N/K)$. Thus
$(\psihat_0,\chi^\star(\psihat_0))$ is an $N$-saddle of $\Gtilde_N$ if this
sign is negative, and an $(N-1)$-saddle otherwise. The same is true for the
index of $\psi=(\psi_1,\dots,\psi_N)$ as a stationary point of~$\Gbar_N$. 

Let $R$ be the so-called\/ \emph{residue}\/ of the periodic orbit of $T$
associated with the stationary point. This residue is equal to
$(2-\Tr(DT^N))/4$, where $DT^N$ is the Jacobian of $T^N$ at the orbit, and
indicates the stability type of the periodic orbit: The orbit is hyperbolic
if $R<0$, elliptic if $0<R<1$, and inverse hyperbolic if $R>1$. It is
known~\cite{MacKayMeiss83} that the residue $R$ is related to the index of
$\psi$ be the identity
\begin{equation}
\label{tmgindex1}
R = -\frac14 \frac{\det(\Hess \Gbar_N(\psi))}{\prod_{j=1}^N (-\sdpar
G{12}(\psi_j,\psi_{j+1}))}\;.
\end{equation}
In our case, $-\sdpar G{12}(\psi_j,\psi_{j+1})$ is always negative, so that
$R$ is positive if $\psi$ is an $(N-1)$-saddle, and negative if  $\psi$ is
an $N$-saddle. 

Now $x^\star(\psihat_0)$ also corresponds to a periodic orbit of the
map~\eqref{tm2}, whose generating function is
$H(x_n,x_{n+1})=\frac12(x_n-x_{n+1})^2 + \frac2\gamma U(x_n)$. The
corresponding  $N$-point generating function is precisely
$(2/\gamma)V_\gamma$.
Since the residue is invariant under area-preserving changes of variables,
we also have 
\begin{equation}
\label{tmgindex2}
R = -\frac1{2\gamma} \frac{\det(\Hess V_\gamma(x^\star))}{\prod_{j=1}^N
(-\sdpar
H{12}(x^\star_j,x^\star_{j+1}))}\;.
\end{equation}
In this case, the denominator is positive. Therefore, $\Hess
V_\gamma(x^\star)$ has an even number of positive eigenvalues if $\psi$ is
an $N$-saddle, and an odd number of positive eigenvalues if $\psi$ is an
$(N-1)$-saddle. 
\end{proof}

We can now complete the proofs of Theorem~\ref{thm_LargeN1} and
Theorem~\ref{thm_LargeN2}. 

\begin{proof}[{\sc Proof of Theorem~\ref{thm_LargeN1}}]
We first recall the following facts, established
in~\cite{BFG06a}. Whenever $\gammat$ crosses a bifurcation value
$\gammat_M$, say from larger to smaller values, the index of the origin
changes from $2M-1$ to $2M+1$. Thus the bifurcation involves a centre
manifold of dimension $2$, with $2M-1$ unstable and $N-2M-1$ stable
directions transversal to the manifold. Within the centre manifold, the
origin repels nearby trajectories, and attracts trajectories starting
sufficiently far away. Therefore, all stationary points lying in the centre
manifold, except the origin, are either sinks or saddles for the reduced
two-dimensional dynamics. For the full dynamics, they are thus saddles of
index $(2M-1)$ or $2M$ (c.f.~\cite[Section~4.3]{BFG06a}), at least for
$\gammat$ close to $\gammat_M$. 

We now return to the twist map in action-angle variables~\eqref{tm16}.  The
frequency $\Omegabar(I)$ being maximal for $I=0$, as $\eps$ increases, new
orbits appear on the line $I=0$, which corresponds to the origin in
$x$-coordinates. Orbits of rotation number $\nu=M/N$ can only exist if
$\eps\Omegabar(0)=\eps\geqs2\pi\nu + \Order{\eps^2}$, which is compatible
with the condition $\gammat<\gammat_M$. 

Consider the case of a winding number $M=1$, that is, of orbits with
rotation number $\nu=1/N$, which are the only orbits existing for
$\gammat_2<\gammat<\gammat_1$. We note that $\gammat>\gammat_2$ implies
$\Delta=2\pi/N\eps > 1/2 - \Order{\eps}$, and thus there exists
$N_1<\infty$
such that the condition $N\geqs N_1$ automatically implies that $\eps$ is
small enough for Corollary~\ref{cor_tmgstat} to hold. Now, 
Proposition~\ref{prop_tmgindex} yields:  
\begin{itemiz}
\item	If $N$ is even, then $K=\gcd(N,2)=2$, and there are $2N$ stationary
points. The points $x^\star(0), x^\star(2\pi/N), \dots$ must be
$2$-saddles,
while the points $x^\star(\pi/N), x^\star(3\pi/N), \dots$ are $1$-saddles; 
\item	If $N$ is odd, then $K=\gcd(N,2)=1$, and there are $4N$ stationary
points. The points $x^\star(0), x^\star(\pi/N), \dots$ must be $1$-saddles,
while the points $x^\star(\pi/2N), x^\star(3\pi/2N), \dots$ are
$2$-saddles. 
\end{itemiz}
Going back to original variables, we obtain the expressions~\eqref{LargeN5}
and~\eqref{largeN7} for the coordinates of these stationary points. The
fact
that they keep the same index as $\gammat$ moves away from $\gammat_M$ is a
consequence of Relation~\eqref{tmgindex2} and the fact that the
corresponding
stationary points of $\smash{\Gtilde_N}$ also keep the same index.  
Finally, Relation~\eqref{largeN6} on the potential difference is a
consequence of Proposition~\ref{prop_tm3}. This proves
Theorem~\ref{thm_LargeN1}. 
\end{proof}

\begin{proof}[{\sc Proof of Theorem~\ref{thm_LargeN2}}]
For larger winding number $M$, one can proceed in an analogous way,
provided $N$ is sufficiently large, as a function of $M$, for the
conditions on $\eps$ to hold. This proves Theorem~\ref{thm_LargeN2}. 
\end{proof}

Finally, Theorem~\ref{thm_stoch2} is proved in an analogous way as
Theorems~2.7 and 2.8 in \cite{BFG06a}, using results from~\cite{FW} (see
also~\cite{Kifer,Sugiura96a}). 


\newpage
\appendix
\section{Jacobi's Elliptic Integrals and Functions}
\label{app_ell}

Fix some $\kappa\in[0,1]$. 
The\/ \emph{incomplete elliptic integrals of the first and second kind}\/ 
are defined, respectively, by\footnote{One should beware of the fact that
some sources use $m=\kappa^2$ as parameter.} 
\begin{equation}
\label{ell1}
\JF(\phi,\kappa) = \int_0^\phi \frac{\6t}{\sqrt{1-\kappa^2\sin^2t}}\;,
\qquad
\JE(\phi,\kappa) = \int_0^\phi \sqrt{1-\kappa^2\sin^2t}\,\6t\;. 
\end{equation}
The\/ \emph{complete elliptic integrals of the first and second kind}\/ 
are given by 
\begin{equation}
\label{ell2}
\JK(\kappa) = \JF(\pi/2,\kappa)\;,
\qquad
\JE(\kappa) = \JE(\pi/2,\kappa)\;.
\end{equation}
Special values include $\JK(0)=\JE(0)=\pi/2$ and $\JE(1)=1$. The integral
of
the first kind $\JK(\kappa)$ diverges logarithmically as
$\kappa\nearrow1$. 

The\/ \emph{Jacobi amplitude}\/ $\am(u,\kappa)$ is the inverse function of
$\JF(\cdot,\kappa)$, i.e.,
\begin{equation}
\label{ell3}
\phi = \am(u,\kappa)
\quad\Leftrightarrow\quad
u = \JF(\phi,\kappa)\;.
\end{equation} 
The three standard Jacobi elliptic functions are then defined as 
\begin{align}
\nonumber
\sn(u,\kappa) &= \sin(\am(u,\kappa))\;,\\
\label{ell4}
\cn(u,\kappa) &= \cos(\am(u,\kappa))\;,\\
\dn(u,\kappa) &= \sqrt{1-\kappa^2\sin^2(\am(u,\kappa))}\;.
\nonumber
\end{align}
Their derivatives are given by  
\begin{align}
\nonumber
\sn'(u,\kappa) &= \cn(u,\kappa)\dn(u,\kappa)\;,\\
\label{ell7}
\cn'(u,\kappa) &= -\sn(u,\kappa)\dn(u,\kappa)\;,\\
\dn'(u,\kappa) &= -\kappa^2\sn(u,\kappa)\cn(u,\kappa)\;.
\nonumber
\end{align}
The function $\sn$ satisfies the periodicity relations
\begin{align}
\nonumber
\sn(u+4\JK(\kappa),\kappa) &= \sn(u,\kappa)\;,\\
\sn(u+2\icx\JK(\sqrt{1-\kappa^2}),\kappa) &= \sn(u,\kappa)\;,
\label{ell7B}
\end{align}
and has simple poles in $u=2n\JK(\kappa)+(2m+1)\icx\JK(\sqrt{1-\kappa^2})$,
$n,m\in\Z$, with residue $(-1)^m/\kappa$. The functions $\cn$ and $\dn$
satisfy similar relations.   Since $\am(u,0)=u$, one has $\sn(u,0)=\sin u$,
$\cn(u,0)=\cos u$ and $\dn(u,0)=1$. As $\kappa$ grows from $0$ to $1$, the
elliptic functions become more and more squarish. This is also apparent
from their Fourier series, given by 
\begin{align}
\nonumber
\frac{2\JK(\kappa)}\pi \sn\biggpar{\frac{2\JK(\kappa)}\pi\psi,\kappa}
&= \frac4\kappa \sum_{p=0}^\infty \frac{\nome^{(2p+1)/2}}{1-\nome^{2p+1}}
\sin\bigpar{(2p+1)\psi} \;, \\
\label{ell5}
\frac{2\JK(\kappa)}\pi \cn\biggpar{\frac{2\JK(\kappa)}\pi\psi,\kappa}
&= \frac4\kappa \sum_{p=0}^\infty \frac{\nome^{(2p+1)/2}}{1+\nome^{2p+1}}
\cos\bigpar{(2p+1)\psi} \;, \\
\nonumber
\frac{2\JK(\kappa)}\pi \dn\biggpar{\frac{2\JK(\kappa)}\pi\psi,\kappa}
&= 1 + 4 \sum_{p=0}^\infty \frac{\nome^p}{1+\nome^{2p}}
\cos\bigpar{2p\psi} \;, 
\end{align}
where $\nome=\nome(\kappa)$ is the\/ \emph{elliptic nome}\/ defined by 
\begin{equation}
\label{ell6}
\nome = \exp\biggset{-\pi\frac{\JK(\sqrt{1-\kappa^2})}{\JK(\kappa)}}\;.
\end{equation}
The elliptic nome has the asymptotic behaviour 
\begin{equation}
\label{ell6a}
\nome(\kappa) = 
\begin{cases}
\vrule height 12pt depth 16pt width 0pt
\dfrac{\kappa^2}{16} + \dfrac{\kappa^4}{32} +
\bigOrder{\kappa^6} 
& \text{for $\kappa\searrow0$\;,} \\
\vrule height 16pt depth 12pt width 0pt
\exp\biggset{\dfrac{\pi^2}{\log\brak{(1-\kappa^2)/16}}}
\biggbrak{1+\biggOrder{\dfrac{1-\kappa^2}{\log^2\brak{(1-\kappa^2)/16}}}}
& \text{for $\kappa\nearrow1$\;.} 
\end{cases}
\end{equation}
We also use the following identities, derived
in~\cite[p.~175]{Jacobi1892}. For $k\geqs1$, 
\begin{equation}
\label{ell8}
\biggpar{\frac{2\JK(\kappa)}\pi}^{2k} 
\sn^{2k}\biggpar{\frac{2\JK(\kappa)}\pi\psi,\kappa}
= \hat c_{2k,0} + \sum_{p=1}^\infty \hat c_{2k,p}
\frac{\nome^{p}}{1-\nome^{2p}}
\cos\bigpar{2p\psi} \;, 
\end{equation}
where the $\hat c_{2k,0}$ are positive constants (independent of $\psi$),
and the other Fourier coefficients are given for the first few $k$ by 
\begin{align}
\nonumber
\hat c_{2,p} &= - \frac4{\kappa^2} (2p)\;, \\
\label{ell9}
\hat c_{4,p} &= \frac4{3!\kappa^4} \biggbrak{(2p)^3 -
4 (2p)(1+\kappa^2)\biggpar{\frac{2\JK(\kappa)}\pi}^2}\;, \\
\hat c_{6,p} &= -\frac4{5!\kappa^6} \biggbrak{(2p)^5 - 20 (2p)^3
(1+\kappa^2)\biggpar{\frac{2\JK(\kappa)}\pi}^2
+ 8 (2p) (8+7\kappa^2+8\kappa^4) \biggpar{\frac{2\JK(\kappa)}\pi}^4}\;.
\nonumber
\end{align}


\section{Proofs of the Fixed-Point Argument}
\label{sec_prtech}


In this appendix, we give the somewhat technical proofs of the fixed-point
argument given in Section~\ref{sec_tmgchi}. We start by proving
Lemma~\ref{lem_tmgff4}, stating a fixed-point equation equivalent to the
stationarity conditions $\tdpar{\Gtilde_N}{\chi_q}=0$, $q\in\cQ$. 

\begin{proof}[{\sc Proof of Lemma~\ref{lem_tmgff4}}]
The definitions~\eqref{tmgff4} of $\alpha_n$ and $\beta_n$ imply, for any
$a, b\geqs 0$, 
\begin{align}
\nonumber
\alpha_n^a &= \sum_{q_1,\dots,q_a\in\cQ} 
\prod_{i=1}^a \rho_{q_i}
\omega^{q_i(n+1/2)} \e^{\icx\psihat_0 q_i/M} \;, \\
\beta_n^b &= \sum_{q'_1,\dots,q'_b\in\cQ} 
\prod_{j=1}^b \frac{-\icx\eps\rho_{q'_j}}{\tan(\pi q'_j/N)} 
\omega^{q'_j(n+1/2)} \e^{\icx\psihat_0 q'_j/M} \;. 
\label{tmgff23:1}
\end{align}
It is more convenient to compute $\tdpar{\Gtilde_N}{\rho_{-q}}$ rather than
$\tdpar{\Gtilde_N}{\chi_q}$. We thus have to compute the derivatives of
$\tilde g_p$ with respect to $\rho_{-q}$ for all $p$. For $p=0$, we have
\begin{equation}
\label{tmgffB5:2}
\dpar{\tilde g_0}{\rho_{-q}} 
= \eps^3 \sum_{n=1}^N G_0'(\Delta+\eps^2\alpha_n)
\dpar{\alpha_n}{\rho_{-q}}\;,
\end{equation} 
where~\eqref{tmgff23:1} shows that $\tdpar{\alpha_n}{\rho_{-q}} =
\omega^{-q(n+ 1/2)}\e^{-\icx\psihat_0 q/M}$. We expand
$G_0'(\Delta+\eps^2\alpha_n)$ into powers of $\eps^2$, and plug
in~\eqref{tmgff23:1} again. In the resulting expression, the sum over $n$
vanishes unless $\sum_i q_i-q$ is a multiple of $N$, say $kN$. This yields 
\begin{equation}
\label{tmgffB5:3}
\dpar{\tilde g_0}{\rho_{-q}} 
= N\eps^3 \sum_{k\in\Z} (-1)^k \e^{\icx k\psihat_0 N/M} 
\sum_{a\geqs0} \frac{\eps^{2a}}{a!} G_0^{(a+1)}(\Delta)
\Gamma^{(a,0)}_{kN+q}(\rho)\;.
\end{equation}
We consider the terms $a=0$ and $a=1$ separately: 
\begin{itemiz}
\item	Since $\smash{\Gamma^{(0,0)}_{kN+q}}(\rho)=\delta_{kN,-q}$ vanishes
for all $k$,
the sum actually starts at $a=1$. 
\item	The fact that
$\Gamma^{(1,0)}_{\ell}(\rho)$ vanishes whenever $\abs{\ell}>N/2$ implies
that only the term $k=0$ contributes, and yields a contribution
proportional
to $-\rho_q$. 
\end{itemiz}
Shifting the summation index $a$ by one unit, we get 
\begin{equation}
\label{tmgffB5:4}
\dpar{\tilde g_0}{\rho_{-q}} 
= N\eps^5 \biggbrak{G_0''(\Delta)\rho_q + 
\sum_{k\in\Z} (-1)^k \e^{\icx k\psihat_0 N/M} 
\sum_{a\geqs1} \frac{\eps^{2a}}{(a+1)!} G_0^{(a+2)}(\Delta)
\Gamma^{(a+1,0)}_{kN+q}(\rho)}\;.
\end{equation}
A similar computation for $p\neq 0$ shows that 
\begin{equation}
\label{tmgffB5:5}
\dpar{\tilde g_p}{\rho_{-q}}\e^{2\icx p\psihat_0} =
N\eps^5  \sum_{k\in\Z} (-1)^k \e^{\icx k\psihat_0 N/M} 
\sum_{a,b\geqs0} \frac{\eps^{2(a+b)}}{a!b!} H^{(a)}_{p,q}(\Delta)p^b
\Gamma^{(a,b)}_{kN+q-2pM}\;.
\end{equation}
Solving the stationarity condition  
\begin{equation}
\label{tmgffB5:1}
0 = \dpar{\Gtilde_N}{\rho_{-q}} 
= \sum_{p=-\infty}^\infty \e^{2\icx p\psihat_0} \dpar{\tilde
g_p}{\rho_{-q}}
\end{equation}
with respect to $\rho_q$, and singling out the term
$a=b=0$ in~\eqref{tmgffB5:5} to give the leading term $\rho^{(0)}$
then yields the result. 
\end{proof}

The following estimates yield sufficient conditions for the operator $\cT$
to be a contraction inside a certain ball, for the norm
$\norm{\cdot}_\lambda$ introduced in~\eqref{tmgff20}. 

\begin{prop}
\label{prop_tmgff2}
There exist numerical constants $c_0, c_1>0$, such that for
any $\lambda<\lambda_0$, and any $N$ such that
$N\e^{-\lambda_0N/2M}\leqs 1/2$, the estimates 
\begin{align}
\label{tmgff21a}
\norm{\cT\rho}_\lambda 
&\leqs \frac{c_1L_0}{\Delta\abs{G_0''(\Delta)}}
\biggbrak{1
+\frac{M}{\Delta^3}\Bigpar{\norm{\rho}_\lambda +
\eps\Delta M \eta(\lambda_0,\lambda)}\eps\norm{\rho}_\lambda}\;, 
\\
\norm{\cT\rho-\cT\rho'}_\lambda 
&\leqs \frac{c_1L_0}{\abs{G_0''(\Delta)}}\frac{M}{\Delta^4}
\Bigbrak{\bigpar{\norm{\rho}_\lambda\vee\norm{\rho'}_\lambda}
+ \eps\Delta M \eta(\lambda_0,\lambda)}\eps
\norm{\rho-\rho'}_\lambda
\label{tmgff21b}
\end{align}
hold with $\eta(\lambda_0,\lambda) =
(\e^\lambda/\lambda_0)\vee(1/(\lambda_0-\lambda))$,  provided $\rho$ and
$\rho'$
satisfy
\begin{equation}
\label{tmgff21c}
\eps \bigpar{\norm{\rho}_\lambda\vee\norm{\rho'}_\lambda}
\leqs c_0 \frac{\Delta^2}M
\biggpar{1\wedge\frac{\lambda_0-\lambda}M\wedge\frac{\lambda}M}\;.
\end{equation}
\end{prop}
\begin{proof}
The lower bound 
\begin{equation}
\label{tmgff21:1}
\frac{\abs{\tan(\pi q/N)}}{\eps} \geqs \frac{\pi\abs{q}}{N\eps}
= \frac{\Delta}{2M}\abs{q}
\end{equation}
directly implies 
\begin{equation}
\label{tmgff21:2}
\abs{H^{(a)}_{p,q}(\Delta)} \leqs L_0 \frac{a!}{r^{a+1}} 
\biggbrak{1+\frac{2M\abs{p}}{\abs{q}}}\e^{-\lambda_0\abs{p}}\;.
\end{equation}
The assumption on $N$ allows $\abs{\rho_q}$ to be bounded by a
geometric series of ratio smaller than $1/2$, which is dominated by the
term
$k=0$, yielding 
\begin{equation}
\label{tmgff21:3}
\norm{\rho^{(0)}}_\lambda 
\leqs \frac{c_2 L_0}{\Delta\abs{G_0''(\Delta)}} 
\e^{-(\lambda_0-\lambda)\abs{q}/2M}
\leqs \frac{c_2L_0}{\Delta\abs{G_0''(\Delta)}}\;,
\end{equation}
where $c_2>0$ is a numerical constant. The fact that
$\Gamma^{(a,b)}_\ell(\rho)$ contains less than $N^{a+b-1}$ terms, together
with~\eqref{tmgff21:1}, implies the bound
\begin{equation}
\label{tmgff21:4}
\abs{\Gamma^{(a,b)}_\ell(\rho)}
\leqs N^{a+b-1} \biggpar{\frac{2M}\Delta}^b
\e^{-\lambda\abs{\ell}/2M}\norm{\rho}_\lambda^{a+b}\;.
\end{equation}
Assuming that $\norm{\rho}_\lambda\leqs c_0\Delta^2/M\eps$ for sufficiently
small $c_0$, it is straightforward to obtain the estimate 
\begin{equation}
\label{tmgff21:5}
\norm{\Phi^{(1)}(\rho,\eps)}_\lambda 
\leqs \frac{c_3L_0}{\abs{G_0''(\Delta)}} \frac{2M}{\Delta^4}
\eps\norm{\rho}_\lambda^2\;. 
\end{equation}
In the sequel, we assume that $q>0$, since by symmetry of the norm under
permutation of $\rho_q$ and $\rho_{-q}$ the same estimates will hold for
$q<0$. The norm of $\Phi^{(2)}(\rho,\eps)$ is more delicate to estimate. We
start by writing 
\begin{equation}
\label{tmgff21:6}
\abs{\Phi^{(2)}_q(\rho,\eps)} \leqs
\frac{L_0}{\abs{G_0''(\Delta)}} \frac1N \sum_{a+b\geqs1}
\bigpar{\eps^2N\norm{\rho}_\lambda}^{a+b} \frac1{r^{a+1}}
\biggpar{\frac{2M}\Delta}^b S_q(b)\;,
\end{equation}
where 
\begin{equation}
\label{tmgff21:7}
S_q(b) = \frac1{b!} \sum_{p\neq0} \abs{p}^{b} \biggpar{1+\frac{2M\abs{p}}q}
\e^{-(\lambda_0-\lambda)\abs{p}}
\sum_{k\in\Z} \exp\biggset{-\frac\lambda{2M}(2M\abs{p}+\abs{kN+q-2Mp})}\;.
\end{equation}
We decompose $S_q(b)=S_q^+(b)+S_q^-(b)$, where $S_q^+(b)$ and $S_q^-(b)$
contain, respectively, the sum over positive and negative $p$. In the
sequel, we shall only treat the term $S_q^+(b)$. The sum over $k$
in~\eqref{tmgff21:7} is dominated by the term for which $kN$ is the closest
possible to $2Mp-q$, and can be bounded by a geometric series. The result
for $p>0$ is
\begin{equation}
\label{tmgff21:7b}
\sum_{k\in\Z} \exp\biggset{-\frac\lambda{2M}(2Mp+\abs{kN+q-2Mp})}
\leqs c_4 \bigpar{\e^{-\lambda p} \wedge \e^{-\lambda q/2M}}\;.
\end{equation}
We now distinguish between two cases. 
\begin{itemiz}
\item	If $q\leqs2M$, we bound the sum over $k$ 
by $\e^{-\lambda p}$, yielding 
\begin{equation}
\label{tmgff21:9}
S^+_q(b) \leqs \frac{c_4}{b!} \frac{4M}q
\sum_{p\geqs1}p^{b+1}\e^{-\lambda_0p}
\leqs 4M c_5 \frac{b+1}{\lambda_0^b}\;.
\end{equation}
Since $\e^{\lambda q/2M}\leqs\e^\lambda$, it follows that 
\begin{equation}
\label{tmgff21:10}
\abs{\Phi^{(2)}_q(\rho,\eps)} \leqs \frac{c_6L_0}{\abs{G_0''(\Delta)}}
\frac{2M^2\e^\lambda}{r^2\Delta\lambda_0} \eps^2\norm{\rho}_\lambda
\e^{-\lambda q/2M}\;.
\end{equation}
\item	If $q>2M$, we split the sum over $p$ at $q/2M$. For $2Mp\leqs q$,
we
bound $(1+2Mp/q)$ by $2$ and the sum over $k$ by
$\e^{-\lambda q/2M}$. For $2Mp>q$, we bound the the sum over $k$ by 
$\e^{-\lambda p}\leqs\e^{-\lambda q/2M}$. This shows 
\begin{equation}
\label{tmgff21:11}
S^+_q(b) \leqs 2c_7M \frac{b+1}{(\lambda_0-\lambda)^b} \e^{-\lambda
q/2M}\;,
\end{equation}
and thus
\begin{equation}
\label{tmgff21:12}
\abs{\Phi^{(2)}_q(\rho,\eps)} \leqs \frac{c_8L_0}{\abs{G_0''(\Delta)}}
\frac{2M^2}{r^2\Delta(\lambda_0-\lambda)} \eps^2\norm{\rho}_\lambda
\e^{-\lambda q/2M}\;.
\end{equation}
\end{itemiz}
Now~\eqref{tmgff21:12} and~\eqref{tmgff21:10}, together
with~\eqref{tmgff21:3}, imply~\eqref{tmgff21a}. The proof
of~\eqref{tmgff21b}
is similar, showing first the estimate 
\begin{equation}
\label{tmgff21:13}
\biggabs{\prod_{i=1}^a \rho_{q_i} - \prod_{i=1}^a \rho'_{q_i}} \leqs
a\bigpar{\norm{\rho}_\lambda\vee\norm{\rho'}_\lambda}^{a-1}
\e^{-\lambda\sum_{i=1}^a\abs{q_i}/2M} \norm{\rho-\rho'}_\lambda
\end{equation}
by induction on $a$, and then 
\begin{equation}
\label{tmgff21:14}
\bigabs{\Gamma^{(a,b)}_\ell(\rho)-\Gamma^{(a,b)}_\ell(\rho')}
\leqs (a+b)
\bigbrak{N(\norm{\rho}_\lambda\vee\norm{\rho'}_\lambda)}^{a+b-1}
\biggpar{\frac{2M}\Delta}^b \e^{-\lambda\abs{\ell}/2M}
\norm{\rho-\rho'}_\lambda\;.
\end{equation}
\end{proof}

It is now easy to complete the proof of Proposition~\ref{cor_tmgff2}. 

\begin{proof}[{\sc Proof of Proposition~\ref{cor_tmgff2}}]
Estimate~\eqref{tmgff21a} for $\norm{\cT\rho}_\lambda$ implies that if 
\begin{equation}
\label{tmgff30:1}
\eps \leqs
\frac{R_0}{\Delta M\eta(\lambda_0,\lambda)}
\wedge \frac{\Delta^3}{2MR_0^2}
\biggpar{\frac{\Delta\abs{G_0''(\Delta)}}{c_1L_0}R_0 - 1}\;,
\end{equation}
then $\cT(\cB_\lambda(0,R_0))\subset\cB_\lambda(0,R_0)$. If in addition 
\begin{equation}
\label{tmgff30:2}
\eps \leqs c_0 \frac{\Delta^2}{MR_0}
\biggpar{1\wedge\frac{\lambda_0-\lambda}M\wedge\frac{\lambda}M}\;,
\end{equation}
then Estimate~\eqref{tmgff21b} for $\norm{\cT\rho-\cT\rho'}_\lambda$
applies for $\rho,\rho'\in\cB_\lambda(0,R_0)$. It is then immediate to
check that $\cT$ is a contracting in $\cB_\lambda(0,R_0)$, as a consequence
of~\eqref{tmgff30:1}. 
Thus the existence of a unique fixed point in that ball follows by Banach's
contraction lemma. Finally, the assertions on the properties of
$\rho^\star$ follow from the facts that they are true for $\rho^{(0)}$,
that they are preserved by $\cT$ and that
$\rho^\star=\lim_{n\to\infty}\cT^n\rho^{(0)}$. 
\end{proof}


\newpage
\small
\bibliography{../BFG}
\bibliographystyle{amsalpha}               

\newpage
\bigskip\bigskip\noindent
{\small 
Nils Berglund \\ 
{\sc CPT--CNRS Luminy} \\
Case 907, 13288~Marseille Cedex 9, France \\
{\it and} \\
{\sc PHYMAT, Universit\'e du Sud Toulon--Var} \\
{\it Present address:} \\
{\sc MAPMO--CNRS, Universit\'e d'Orl\'eans} \\
B\^atiment de Math\'ematiques, Rue de Chartres \\
B.P. 6759, 45067 Orl\'eans Cedex 2, France \\
{\it E-mail address: }{\tt berglund@cpt.univ-mrs.fr}

\bigskip\noindent
Bastien Fernandez \\ 
{\sc CPT--CNRS Luminy} \\
Case 907, 13288~Marseille Cedex 9, France \\
{\it E-mail address: }{\tt fernandez@cpt.univ-mrs.fr}

\bigskip\noindent
Barbara Gentz \\ 
{\sc Weierstra\ss\ Institute for Applied Analysis and Stochastics} \\
Mohrenstra{\ss}e~39, 10117~Berlin, Germany \\
{\it Present address:}\\
{\sc Faculty of Mathematics, University of Bielefeld} \\
P.O. Box 10 01 31, 33501~Bielefeld, Germany \\
{\it E-mail address: }{\tt gentz@math.uni-bielefeld.de}

}


\end{document}